\documentclass[11pt]{article}
\usepackage[utf8]{inputenc}
\usepackage{tipa}
\usepackage{dsfont}
\usepackage{amssymb}
\usepackage{stmaryrd}
\usepackage{dirtytalk}
\usepackage{amsmath,amscd}      
\usepackage{amsfonts}       
\usepackage{amsthm}         
\usepackage{bbding}         
\usepackage{bm}             
\usepackage{graphicx}       
\usepackage{fancyvrb}       
\usepackage{indentfirst}    
\usepackage{natbib}         
\usepackage[nottoc]{tocbibind} 
\usepackage{icomma}         
\usepackage{dcolumn}        
\usepackage{booktabs}       
\usepackage{paralist}       
\usepackage{xcolor}         
\usepackage{ dsfont }
\usepackage{blkarray}
\usepackage{lmodern}
\usepackage[all,cmtip]{xy}
\usepackage[T1]{fontenc}
\usepackage{textcomp}
\theoremstyle{plain}
\newtheorem{veta}{Věta}

\newtheorem{Lemma}[veta]{Lemma}

\newtheorem{Thm}[veta]{Theorem}
\newtheorem{Prop}[veta]{Proposition}
\newtheorem{Ex}[veta]{Example}

\newtheorem{Cor}[veta]{Corollary}

\newtheorem{Def}[veta]{Definition}

\theoremstyle{remark}

\newtheorem*{remark}{Remark}

\theoremstyle{plain}

\newenvironment{Proof}{
  \par\medskip\noindent
  \textit{Proof}.
}{
\newline
\rightline{$\qedsymbol$}
}

\usepackage{tikz-cd}

\newcommand{\m}{\mathfrak{m}}

\usepackage[utf8]{inputenc}
\usepackage{tikz}
\usetikzlibrary{positioning, calc, arrows.meta}
\usepackage{blindtext}
\title{Semiperfect rings with a Nakayama permutation: A survey of Double annihilator property and Size condition}
\author{Dominik Krasula
}

\begin{document}

\maketitle


\textbf{Abstract} 

For a semiperfect ring with essential socles, the Double annihilator property encodes that the top and socle have anti-isomorphic lattices of submodules, whereas the Size condition encodes that they are isomorphic as modules. Interest in both concepts, particularly for finite rings, was revived by coding theory, where they characterise QF rings and Frobenius rings, respectively. However, their shared origins date back to the work of T. Nakayama.

We study these concepts through the lens of the Nakayama permutation, an invariant initially used to define (quasi-)Frobenius rings. We propose semiperfect rings as the setting for this study, treating them as the natural generalisation of finite rings, because they possess the characteristic decomposition of unity preserved by projection onto a semisimple top. This allows us to extend the utility of the Nakayama permutation beyond the classical Artinian setting.

By analysing the Nakayama permutation in this broader context, we show that many classical properties of (quasi-)Frobenius rings are not exclusive to the finite case, but are special cases of the general behaviour of semiperfect rings with essential socles.

We illustrate these results using B. J. Müller's representation of semiperfect rings as rings of formal matrices. The clear description of socles and tops in this setting provides a straightforward method for constructing counterexamples, such as quasi-Frobenius rings that are not Frobenius.

\smallskip

\textbf{Keywords:} semiperfect rings with essential socles; annihilator of an ideal; Nakayama permutation; Kasch rings

\smallskip

\textbf{Mathematics Subject Classification} 16L30; 16L60; Secondary: 16S50

\smallskip

\noindent \textbf{Author:} RNDr. Dominik Krasula

Charles University, Faculty of Mathematics and Physics, 

Department of Algebra

Sokolovská 49/83, 186 75 Praha 8, Czech Republic

krasula@karlin.mff.cuni.cz

ORCID: 0000-0002-1021-7364

\smallskip 

\noindent This work is a part of project SVV-2025-260837.

\noindent This research was supported by the grant GA ČR 23-05148S from the Czech
Science Foundation.

\section{Introduction}

A ring is called \textit{Frobenius} if it is an artinian ring such that $soc(R_R)\cong top(R_R)$ and $soc(_RR)\cong top(_RR)$. For a finite ring, it is enough to assume one of the isomorphisms (Thm. \ref{ThmIovCyclic}). A ring is \textit{quasi-Frobenius}, or \textit{QF} for short, if it is Morita equivalent to a Frobenius ring. 

These two concepts were introduced by T. Nakayama in [Nak. 41a, Chap. II] by means of a \textit{Nakayama permutation} (Def. \ref{DefNak}) generalising concepts previously defined for fin. dim. algebras [Nak. 39]. This was motivated by the observation that QF rings are characterised by the  \textit{Double annihilator property} (Sec. \ref{SecDouble}). This property provides an anti-isomorphism between the lattices of submodules of the ring's top and socle (Sec. \ref{SecDIsom}). Conversely, the existence of any such anti-isomorphism implies the Double annihilator property [Nak. 41b]. A stronger version of this property, a \textit{cardinality condition}, characterises Frobenius rings (Sec. \ref{SecCard}).

In subsequent decades, QF rings were extensively studied, particularly for their homological properties. They are characterised as artinian self-injective rings. A ring is QF if and only if all right (and hence also left) projective modules are injective [Fai. 66, Thm. A], or equivalently, if right (and hence also left) injective modules are projective [FW 67, Thm. A]. This line of study revealed a strong connection between structural and homological properties and notably led to the celebrated \textit{Faith-Walker Theorem} [FW 67, Thm. B]. An extensive survey of the classical literature can be found in [Fai. 76, Notes on Chap. 24].

\smallskip 

In [Wood 99], J. A. Wood showed that both MacWilliams theorems hold for the class of finite Frobenius rings. S. T. Dougherty writes in his monograph [Dou. 17] that \say{\textit{Generally, when studying codes over rings, a blanket assumption is made that all
rings serving as alphabets for codes are finite Frobenius rings.}} Other researchers proposed \textit{quasi-Frobenius bimodules} as alphabets [GNW 04], but these are still tied to results on QF rings. 

The newly established connection with coding theory has led to a revived interest in the theoretical properties of (finite) QF rings, specifically the Double annihilator property and the Size condition. However, the extensive theoretical literature on QF rings is not fully leveraged in addressing practical problems in coding.

This survey revisits classical results on QF rings from the perspective of the aforementioned properties used in coding theory, interpreting them through the lens of the Nakayama permutation. We propose \textit{semiperfect rings with a Nakayama permutation (and an essential right socle)} as the setting for this study. Throughout the text, we demonstrate several properties of finite (quasi-)Frobenius rings as straightforward special cases of the results on semiperfect rings with a Nakayama permutation. This might be somewhat surprising, since in the study of MacWilliams conditions, treating the finite case separately proved necessary [Iov. 22, Cor. 0.3].

\smallskip 

We review the theory of semiperfect rings used throughout the text in Sec. \ref{SecSemp}. We then study the definition of the Nakayama permutation (Sec. \ref{SecDef}), its connection to the Double annihilator property (Sec. \ref{SecDouble}), and the Size condition (Sec. \ref{SecCard}). We accompany the text with examples constructed using a specific version of formal matrix rings, as defined in Sec. \ref{Secsempmatrix}.

\subsection{Semiperfect rings and a Nakayama permutation}

In a finite ring $R$, one can express its unity as a sum of pairwise orthogonal primitive idempotents, and for any primitive idempotent, the corner ring $eRe$ is local.  B. J. Müller showed that this property characterises semiperfect rings [Mül. 70, Thm. 1]. This, in particular, implies that the concept of \textit{a checkered matrix ring} used for finite rings [McD. 74, Sec. X.] can be generalised to semiperfect rings (Sec. \ref{Secsempmatrix}). 

We detail the structural theory of semiperfect rings and interpret it as a direct generalisation of standard results on finite rings, as presented in B. R. McDonald's classical monograph [McD. 74] (Sec. \ref{SecSemp}), though some notational conventions differ. 

The decomposition of the unity establishes semiperfect rings as a natural setting to define a Nakayama permutation. As we demonstrate throughout this survey, many classical results on QF rings can be immediately generalised to semiperfect rings with essential socles.  Nakayama could not formulate such generalisations as the concept of a semiperfect ring was yet to be established.

\smallskip

H. Bass obtained semiperfect and (right) perfect rings as natural generalisations of semiprimary rings [Bass 60, Thms. P, 2.1]. Projective and injective modules over perfect rings share most of their fundamental properties with those over finite rings; indeed, a right perfect self-injective ring is QF. This follows from results by B. Osofsky and T. Kato, as observed by C. Faith. It remains an open problem whether it is sufficient to assume self-injectivity only on one side; see the survey [FV 02]. Faith originally formulated it as an exercise to prove the negative [Fai. 76, Ex. 22.16]. 

Semiprimary rings are, in turn, generalisations of (right) artinian rings, which are a common generalisation of both fin. dim. algebras over a field and finite rings. E. Artin originally studied rings with both maximum condition (now called \textit{right Noetherian}) and minimum condition. His seminal article [Art. 27] is often attributed to proving that rings with those two conditions are \textit{semiprimary}, i.e., have a nilpotent Jacobson radical [Art. 27, Satz 4], and are semilocal [Art. 27, Satz 11]. In this view, perfect rings generalise them by requiring that the radical is T-nilpotent.

C. Hopkins and J. Levitzki proved that a maximum condition is unnecessary for establishing that the radical is nilpotent [Hop. 39, 1.4], [Lev. 40, Thm. 6] and that the top is semisimple [Lev. 40, Thm. 10]. This showed that right artinian rings are semiprimary [Hop. 38, 3.2]. It is well known that a Noetherian semiprimary ring is artinian [Hop. 39, 6.7], cf. [Fai. 76, Prop. 18.12]. A right perfect left Noetherian ring is left artinian [Lam 99, Thm. 23.20(6]. Further conditions, where the Noetherian condition implies artinian, can be found in Sec. \ref{SecDefNoether}.

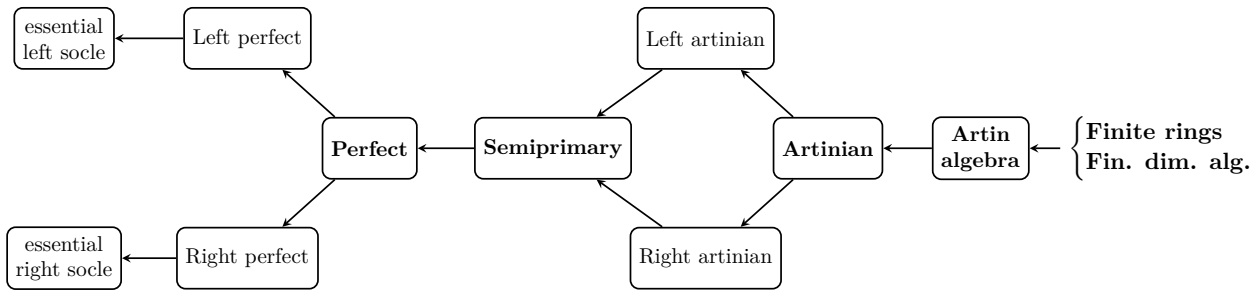
\begin{figure}[ht]
    \centering
    \resizebox{\textwidth}{!}{
    \begin{tikzpicture}[
        >=stealth,
        block/.style={
            draw=black, 
            rectangle, 
            rounded corners, 
            thick, 
            align=center, 
            fill=white, 
            font=\small,
            inner sep=4pt,
            minimum height=1cm
        },
        arrow/.style={->, thick, draw=black}
    ]

    \node[block] (ArtinAlg) at (7.5, 0) {\textbf{Artin} \\\textbf{algebra}};
    
    \node[align=left, anchor=west] (ArtinContent) at (8.8, 0) {
        $\left\{\begin{aligned}
            &\textbf{Finite rings} \\[-0.5ex]
            &\textbf{Fin. dim. alg.}
        \end{aligned}\right.$
    };
    
    \draw[arrow] (ArtinContent.west) -- (ArtinAlg.east);

    \node[block] (Artinian) at (5.0, 0) {\textbf{Artinian}};

    \node[block] (LArt) at (3.0, 1.8) {Left artinian};
    \node[block] (RArt) at (3.0, -1.8) {Right artinian};

    \node[block] (SemiPrim) at (0.5, 0) {\textbf{Semiprimary}};

    \node[block] (Perfect) at (-2.5, 0) {\textbf{Perfect}};

    \node[block] (LPerf) at (-4.5, 1.8) {Left perfect};
    \node[block] (RPerf) at (-4.5, -1.8) {Right perfect};

    \node[block] (LSemiPerf) at (-7.5, 1.8) {essential \\left socle};
    \node[block] (RSemiPerf) at (-7.5, -1.8) {essential\\right socle};


    \draw[arrow] (ArtinAlg) -- (Artinian);

    \draw[arrow] (Artinian) -- (LArt);
    \draw[arrow] (Artinian) -- (RArt);

    \draw[arrow] (LArt) -- (SemiPrim);
    \draw[arrow] (RArt) -- (SemiPrim);

    \draw[arrow] (SemiPrim) -- (Perfect);

    \draw[arrow] (Perfect) -- (LPerf);
    \draw[arrow] (Perfect) -- (RPerf);

    \draw[arrow] (LPerf) -- (LSemiPerf);
    \draw[arrow] (RPerf) -- (RSemiPerf);

    \end{tikzpicture}
    }
    \label{FigRingImplications}
\end{figure}
We are particularly interested in an intermediate class of rings, namely \textit{semiperfect rings with an essential right socle}.  This means that all right ideals contain a simple right ideal. A left perfect ring can be characterised as a semiperfect ring such that all right \textit{modules} have a simple submodule [Fai. 76, Thm. 22.29(3)].  The assumption that a ring is perfect is unnecessarily strong for our purposes, as we only deal with the intrinsic properties of a ring.  

When possible, we relax the assumption that the socle is essential. We often encounter classes of rings where the essentiality of (one or both) socles implies that they coincide (Thm. \ref{PropQFKaschEss}). Several key results hold for any semiperfect ring with coinciding socles (Thm. \ref{ThmAnn1}, \ref{ThmAnnMain}, Prop. \ref{PropIsomfields}, Lemma \ref{LemmaAnnofsimple}). Nevertheless, the essentiality assumption arises naturally, as we are primarily interested in ring properties derived from the properties of the socle.

\subsection{Main results}

If a semiperfect ring $R$ has a Nakayama permutation, it follows immediately that all ideals generated by idempotents have simple socles, i.e., it is a \textit{QF-2 ring}, and that it is \textit{Kasch}. We show that a ring $R$ with an essential right socle has a Nakayama permutation if and only if it is a Kasch QF-2 ring (Thm. \ref{PropQFKaschEss}) and its socles then coincide.

Semiperfect rings with a Nakayama permutation and coinciding socles can be characterised as Kasch minannihilator rings. The Kasch condition can also be viewed as an annihilator condition  (Thm. \ref{ThmKasch}) and semiperfect Kasch minannihilator rings are characterised by the fact that annihilators provide duality between lattices of submodules of their socles and top (Thm. \ref{ThmAnn1}). If both socles are essential, any such duality is given by the annihilator maps (Thm. \ref{ThmAnti-isom}). 

 Thm. \ref{ThmAnnMain} shows that \textit{Generalised dimension condition} (a version of Wood's Size condition) on semisimple ideals characterises when a Nakayama permutation respects multiplicities (a generalisation of Frobenius rings). For rings with essential socles, this combinatorial condition is equivalent to the classical structural property that $soc(R)\cong top(R)$ (Cor. \ref{CorCardmain}), extending results by Nakayama, Wood and T. Honold.

We view semiperfect rings with a Nakayama permutation respecting multiplicities as the natural generalisation of Frobenius rings, while those with an arbitrary Nakayama permutation generalise QF rings. Indeed, the latter class consists precisely of rings Morita equivalent to the former. As Morita equivalence is characterised by varying the multiplicities of isomorphism classes of primitive ideals in the decomposition of the ring (Sec. \ref{Secsempbasic}), this perspective clarifies why the module isomorphism $top(R) \cong soc(R)$ is characterised by the stricter condition, whereas the lattice anti-isomorphism survives in the general case. Varying the multiplicities disrupts the module isomorphism because the simple composition factors of the socle whose multiplicities change ($V_{\pi(k)}$) are not isomorphic to those varying in the top ($V_k$). However, the lattice anti-isomorphism persists. This is because the simple modules $V_k$ and $V_{\pi(k)}$ belong to the same cycle of the Nakayama permutation and thus share isomorphic endomorphism rings (Prop. \ref{PropIsomfields}). Consequently, the underlying geometries of the corresponding homogeneous components remain compatible, preserving the duality of the full lattices even when the modules themselves are no longer isomorphic.

\subsection{Notation}

Whenever we specify an $(R\text-S)$-bimodule structure on some abelian group $B$, viewing $B$ \textit{as a right (or left) module} implies the restriction of this bimodule structure (i.e., it can be formalised by a forgetful functor $R\text- Mod\text- S\to  Mod\text-R$).

When we refer to a \textit{side-specific property} of a bimodule (or ring) without specifying the side, we always mean \textit{on both sides.}. For example, by an \textit{artinian ring}, we mean that it satisfies both the left and right minimum conditions

Whenever we refer to \say{\textit{a semiperfect ring $R$}}, we always assume the notation introduced in Sections \ref{SecSempNot} and \ref{Secsemptop}. In particular, $S_r$ is its right socle, $e_1, \dots, e_n$ is a fixed set of basic idempotents, and so on.

\section{Semiperfect rings and formal matrix representations}\label{SecSemp}

This section presents the general theory of semiperfect rings used throughout the text and outlines the notational conventions of this article. We discuss how to construct examples of QF rings that are not Frobenius (Sec. \ref{secsempex}). 

While the Nakayama permutation is not properly defined until Definition \ref{DefNak}, we formulate the key Lemma \ref{LemmaSocleStruct} that allows us to compare the lattice structure of the socle and the top of a QF ring.

\smallskip 

We start with the finite decomposition of unity, generalising properties of finite rings (Sec. \ref{SecSempNot}). 

Given the article's primary focus on the relationship between the top and the socle, we describe them in detail. Semiperfect rings can be characterised by the properties of their tops (Sec.  \ref{Secsemptop}). The properties of the socle are less restrictive, yet it remains well-behaved with respect to the decomposition of unity (\ref{secsempsocle}). 

We then proceed to describe Müller's formal matrix ring representation of semiperfect rings. Over finite rings, this coincides with what McDonald calls \textit{checkered matrix rings}. We describe the construction (Sec. \ref{Secsempmatrix}) and discuss how to reduce the study to basic rings (Sec. \ref{Secsempbasic}). The construction is then illustrated through a series of examples (Sec. \ref{secsempex}), and we provide detailed descriptions of the representations of simple modules (Sec. \ref{SecsempSimples}).  

The section ends with a survey on Morita duality (Sec. \ref{SecDual}). While the concept of Morita duality is not used in our proofs, we include a brief discussion to provide contextual background for these results.

\subsection{Basics and notation of semiperfect rings}\label{SecSempNot}

This section shows how the decomposition of the ring's unity as a sum of local idempotents induces a decomposition of the ring and provides a complete characterisation of finitely generated (f.g) projective modules, cf. [McD. 74, Thm. XII.4]. In particular, we define invariants $m$ and $n$, a basic set of idempotents $e_1, \dots, e_n$ with multiplicities $\mu_1, \dots, \mu_n$ and index sets $I_1, \dots, I_n$,  and primitive and basic right ideals.

We use [Lam 91, Chapters 7, 8] as the reference for the properties of idempotents in semiperfect rings, noting that these properties largely coincide with those in [McD. 74, Sec. VII].

\smallskip 

An idempotent $e\in R$ is called a \textit{local idempotent} if $End_R(eR)$ is \textit{local}, i.e., it has only one right maximal ideal. This ring is isomorphic to the \textit{corner ring} $eRe$. Müller observed that a ring $R$ is \textit{semiperfect} if and only if its unity $1_R$ can be written as a sum of finitely many pairwise orthogonal local idempotents [Mül. 70, Thm. 1]. In a semiperfect ring, an idempotent is local if and only if it is \textit{primitive}.

This decomposition of $1_R$ is unique up to conjugation and reordering. The total number of these idempotents, $m$, is a fixed invariant of $R$. In the context of this article, we shall call this $m$ the \textit{order of $R$} drawing a parallel with the \textit{order of a matrix}. Note that this is the same as \textit{Krull-Schmidt length} of $R_R$.

Let $R$ be a semiperfect ring of order $m$, and let $V_1,\dots, V_n$ be a list of all simple right $R$-modules up to isomorphism, where $n\leq m$. In Section \ref{SecsempSimples}, we describe a canonical representation of simple right modules.  This $n$ is sometimes referred to as the \textit{height of $R$}; however, the term \textit{height of a ring} is not standardised and is used in the literature with various meanings.

We can order a complete set of pairwise orthogonal local idempotents $e_1,\dots, e_m$ in $R$, such that for each $k\leq n$, the principal right ideal $e_kR$ is isomorphic to the projective cover of $V_i$. Then $e_1R,\dots, e_nR$ is a complete set of f.g. indecomposable projective right $R$-modules up to isomorphism, and we will call them the \textit{basic right ideals}. Because their rings of endomorphisms are local, any f.g. projective right module has a unique decomposition as a direct sum of direct powers of basic right ideals [AF 92, Thm. 12.6].

A right ideal generated by any local idempotent is called a \textit{primitive right ideal} (or \textit{indecomposable principal module}). We reserve the term \textit{basic right ideals} specifically for the chosen representatives $e_1R,\dots, e_nR$. 

The idempotents $e_1,\dots, e_n$ are called a \textit{basic set of idempotents}. We associate with this set \textit{multiplicities}, a set of positive integers $\mu_1,\dots, \mu_n$ such that 
 \begin{gather*}
     R_R\cong (e_1R)^{\mu_1}\oplus \dots \oplus  (e_nR)^{\mu_n}  \\ 
          _R R\cong (Re_1)^{\mu_1}\oplus \dots \oplus  (Re_n)^{\mu_n},
 \end{gather*}
in the category of right and left modules, respectively.

Furthermore, to each $k\leq n$, we assign the set $I_k\subseteq \{1, 2, \dots, m\}$ consisting of indices such that $e_kR\cong e_jR$ if and only if $j\in I_k$. In particular $\mu_k=\mid I_k\mid$. Note that this is equivalent to $Re_k\cong Re_j$.

\subsection{Top of the ring}\label{Secsemptop}

A ring $R$ is semiperfect if and only if it is semilocal and idempotents lift modulo the Jacobson radical $J(R)$.  We use this characterisation to provide a detailed description of the right ideals containing $J(R)$. We also fix the notation for the \textit{skew fields} (also called \textit{division rings}) $\m_1, \dots, \m_n$,  simple right (left) modules $V_k$ ($V'_k$), and lattices of submodules $\mathcal{L}(M)$.

A ring $R$ is \textit{semilocal} if $top(R):=R/J(R)$ is a semisimple ring.  A ring with finitely many maximal ideals is semilocal; however, as we will see, the converse is generally not true in the noncommutative case. Note that the terminology here differs from McDonald's \textit{semi-local rings} [McD 74, Thm. XII.16]. But the Wedderburn-Artin theory used to study the top coincides with [McD. 74, Thms. VIII.3,4].

\smallskip

For a ring $R$ and a right (left) $R$-module $M$ ($M'$), we denote \textit{lattice of submodules} of $M$ ($M'$) by $\mathcal{L}_R(M)$ ($_R\mathcal{L}(M')$). When clear from the context, we omit the subscript $_R$. 

Recall that for any submodule $N\subseteq M$, the lattice $\mathcal{L}(M/N)$ is isomorphic to the interval $[N,M]$ in ${\mathcal{L}(M)}$ via the natural projection  mapping $N\subseteq N'\subseteq M$ to $N'/N$.

In particular, we have a lattice isomorphism \[
\begin{array}{rcl}
    [J(R), R]_{\mathcal{L}_R(R)} & \longleftrightarrow & \mathcal{L}_R(top(R)) \\[1ex]
    I & \longmapsto & I/J(R) \\[1ex]
    I'+J(R) & \mapsfrom & I'
\end{array}
\]
and similarly for the left submodules. Furthermore, any $R$-submodule of $top(R)$ is a $top(R)$-module, and vice versa. 

The structure of $top(R)$ as a ring is well-known by means of the Wedderburn–Artin theorem. There exist skew fields $\m_1, \dots, \m_n$ such that there exists a ring isomorphism \[top(R)\cong M_{\mu_1}(\m_1)\times \dots \times  M_{\mu_n}(\m_n), \]
where $\m_k\cong top(e_kRe_k)$ is isomorphic to the endomorphism ring of a simple right $R$-module $V_k\cong top(e_kR)$.

Note that this gives a natural bijection mapping (an isomorphism class of) a simple right module $V_k$ to a simple left module $V'_k$.

\smallskip 

Because idempotents lift modulo the radical, we can strengthen the above correspondence to describe all ideals containing the radical.
\begin{Lemma}\label{LemmaIdealsAboveJ}
Let $R$ be a semiperfect ring.

Then any right ideal $I$ containing $J(R)$ can be written as $(1-f)R+J(R)$ for some idempotent $f\in R$.    

Furthermore, $I$ is maximal if and only if $f$ is primitive. 
\end{Lemma}
\begin{Proof}
Since the ring $top(R)$ is semisimple, it is von Neumann regular [Lam 91, Cor. 24], and thus all its right ideals are generated by an idempotent [Lam 91, Thm. 23].  

    The right ideal $I/J(R)$ is generated by an idempotent $\bar{e}top(R)$, so $I/J(R)=\bar{e}top(R)$. Since idempotents lift modulo $J(R)$, there exists an idempotent $e\in R$ such that $I=J(R)+eR$. We then substitute $f=1-e$.

The final part follows from the corresponding fact about ideals in the ring $top(R)$.   
\end{Proof}

For a semisimple module, the lattice of submodules can be fully described in terms of its simple submodules. The following result appears throughout the literature, see [Kra. 24, Sec. 3.1]
\begin{Lemma}\label{LemmaUrelated}
    Consider a family of semisimple modules $M_i$.

    Then the lattice of submodules $\mathcal{L}(\oplus M_i)$ is isomorphic to $\prod\mathcal{L}(M_i)$ if and only if no pair of distinct modules $M_i$ and $M_j$ shares an isomorphic simple submodule.
\end{Lemma}

In particular, for a semiperfect ring $R$, we have
\begin{gather*}
    \mathcal{L}_R(top(R))\cong \prod_{k=1}^n \mathcal{L}(V_k^{\mu_k}).
\end{gather*}

\subsection{Properties of the socle} \label{secsempsocle}

This section describes the properties of socles, focusing on properties specific to semiperfect rings, namely the relation between the $soc(R_R)$ and the socles of primitive right ideals. We fix the notation for the right and left socle $S_r$ and $S_l$ and their homogeneous components $S_k$ and $S'_k$. We end by describing the lattice of submodules of the socle of a ring with a Nakayama permutation (Lemma \ref{LemmaSocleStruct}).

\smallskip

We denote the right socle by $S_r:=soc(R_R)$ and the left socle by $S_l:=soc(_RR)$. These are two-sided ideals. By \textit{essential right socle}, we mean that $S_r$ is an \textit{essential right ideal}, i.e., it has a nonzero intersection with all nonzero right ideals. For rings appearing in most of the main results, $S_l$ and $S_r$ coincide; in such cases, we will use the notation $soc(R)$.

To each simple right module $V_k$, we assign a \textit{right homogeneous component $S_k$} of the right socle $S_r$, which is the sum of all simple right ideals isomorphic to $V_k$. If the ring is not right Kasch, some homogeneous components may be zero. We denote left homogeneous components by $S'_k$.

\begin{Lemma}
    Let $R$ be a ring and let $S$ denote a homogeneous component of its right socle.

    Then the right ideal $S$ is a left ideal. 
\end{Lemma}
\begin{Proof}
     For any nonzero $t\in S$, the right module $tR\subseteq S$ is a direct sum of isomorphic copies of some fixed simple right module, and so is its image under any nonzero homomorphism. Now, for any $r\in R$, the left multiplication by $r$ is a homomorphism of right modules. Hence,  $rtR$ is contained in $S$.
\end{Proof}

\begin{Lemma}\label{LemmaSoceR}
    Let $R$ be a semiperfect ring. 

    (1) If $T$ is a simple right ideal, then it is isomorphic to a simple subideal of $e_iR$ for some $i\leq m$.

    (2) If $f\in R$ is an idempotent, then  $fS_r=fR\cap S_r=soc(fR)$.

    (3) $S_r=\oplus^m_{i=1} soc(e_iR)$.
\end{Lemma}
We will prove that in rings with a Nakayama permutation and $S_r=S_l$, all right semisimple ideals are of the form $fsoc(R)$ for some idempotent $f$ (Cor. \ref{CorDmain})
\begin{Proof}
(1) Because $1_RT=T$, there exists $i\leq m$ such that $e_iT\neq 0$. Now define a homomorphism of right modules $T\to e_iT$ by $t\mapsto e_it$. It is a nonzero homomorphism of a simple module; thus, it is injective and hence $e_iR$ contains an isomorphic copy of $T$.     

(2) Let $r\in R$ such that $fr\in S_r$, then $fr=f^2r$ is in $fS_r$. Now take $s\in S_r$, because $S_r$ is an ideal, $fs\in S_r$, and the second inclusion follows.

    The second equality can be proved without assuming that $f$ is idempotent, as both sides represent the largest semisimple submodule of $fR$.

(3) Clearly, socles of right primitive ideals are contained in $S_r$. For the opposite direction, we repeat the argument from (1). Let $T$ be a semisimple right ideal, then $T=e_1T+\dots e_mT$. Now each $e_iT$ is either zero or a semisimple submodule of $e_iR$.  
\end{Proof}

We now discuss the lattice-theoretic properties of general semisimple submodules. We use [Bir. 48] as a reference.  
\begin{Lemma}\label{LemmaBir} Let $R$ be a semiperfect ring and $M$ a semisimple right module of finite composition length.

(1) $\mathcal{L}_R(M)$ is a modular lattice of finite length. 

(2) If $M$ is isomorphic to a direct product of copies of $V_k$ for some $k\leq n$, then $\mathcal{L}_R(M)=\mathcal{L}_{M_{\mu_k}(\m_k)}(M)$ is a projective geometry and is thus indecomposable as a lattice.

(3) The decomposition of $M$ as a product of its homogeneous components induces a decomposition of $\mathcal{L}_R(M)$ as a product of indecomposable lattices.

(4) The rank of  $\mathcal{L}_R(M)$ equals the composition length of $M$.
\end{Lemma}
\begin{Proof}
    (1) All lattices of submodules are modular [Bir. 48, Thm. V.1], and a composition length of a module corresponds to the length of the lattice of submodules.

     (2) Let $I$ be the maximal (two-sided) ideal of $R$ such that $R/I$ corresponds to $M_{\mu_k}(\m_k)$ in a decomposition of $top(R)$ as a ring. Then $I$ annihilates all simple modules isomorphic to $V_k$ and thus any finite direct sum of isomorphic copies of $V_k$ is a free $R/I\cong M_{\mu_k}(\m_k)$-module and hence a \textit{projective geometry} in the sense of [Bir. 48, Ex. VII.2] and thus indecomposable as a lattice [Bir. 48, Cor. to Thm. VIII.6]

     (3) This follows from Lemma \ref{LemmaUrelated}.

     (4) By (3) and uniqueness of such decompositions [Bir. 48, Thm. VI. 12] we can assume that $M$ is a direct product of isomorphic copies of some fixed simple module $T$. 
     
     The \textit{rank} is the maximum number of independent \textit{atoms} (also called \textit{points} in [Bir. 48]). Atoms of $\mathcal{L}_R(M)$ correspond to simple submodules, and the conclusion then follows from [Bir. 48, Lemma VII. 1]
\end{Proof}
In a ring with a Nakayama permutation, primitive ideals are up to isomorphism determined by their simple socles. Lemma \ref{LemmaUrelated} then gives:
\begin{Lemma}\label{LemmaSocleStruct}
Let $R$ be a semiperfect with a Nakayama permutation $\pi$.

Then $\mathcal{L}_R(S_r)\cong \prod_{k=1}^n \mathcal{L}((V_{\pi(k)})^{\mu_{k}})$. 
\end{Lemma}

\subsection{Formal matrix rings}\label{Secsempmatrix}

This section describes Müller's formal matrix ring representation of semiperfect rings. We also define a special case, the \textit{trivial matrix rings}, which will be used for the construction of examples throughout the text. Müller's construction is a natural generalisation of checkered matrix rings [McD. 74, Sec. X], but McDonald's construction assumes a different ordering of a complete set of primitive idempotents than we do.

Formal matrix rings have appeared repeatedly in the literature. Classic examples are triangular matrix rings. A prominent example by S. U. Chase [Cha. 61, Prop. 3.1] provides a construction of a right semihereditary ring that is not left semihereditary. E. L. Green described the module theory over such rings in [Gre. 82], which we recall in Sec. \ref{SecsempSimples}.  A comprehensive study of the general properties of formal matrix rings can be found in [KT 17].

\smallskip

Let $R$ be a semiperfect ring. Then we have the following series of isomorphisms:
\[R\cong End_R(R_R)=Hom(\bigoplus_{1 \leq i \leq m } e_iR, \bigoplus_{1 \leq i \leq m } e_iR)\cong \bigoplus_{1\leq i, j \leq m} Hom_R(e_jR, e_iR)\cong  \bigoplus_{1\leq i, j \leq m} e_iRe_j \]
in the category of abelian groups, sometimes referred to as a \textit{Peirce decomposition of $R$}.

Müller observed in [Mül. 70, Sec. 2] that this allows us to represent $R$ as a set of square matrices of order $m$, such that for $i, j\leq m$, an element in the $(i,j)$-th entry is an element of the $R_i\text-R_j$-bimodule $e_iRe_j$, where $R_i:=e_iRe_i$ is a local ring. The standard matrix addition and multiplication correspond to the addition and multiplication in $R$. This representation is not unique, as the decomposition of $1_R$ is determined only up to reordering and conjugation by a unit.

In the opposite direction, let $m$ be a positive integer and $R_1,\dots, R_m$ a set of local rings, and for $i, j, k\leq m$, let $B_{i,j}$ be an $(R_i\text-R_j)$-bimodule with $B_{i,i}:=R_i$. Then we can define a \textit{formal matrix ring}, also called \textit{a ring of generalised matrices}, where each element is an $m\times m$ table with an element of the \textit{coordinate bimodule} $B_{i,j}$ in the $(i, j)$-position, using the standard matrix addition and multiplication. 
\[
R:=\begin{blockarray}{ccccc} 
    Re_1 & Re_2 & \cdots & Re_m \\[1ex]
    \begin{block}{(cccc)c}
        R_1 & B_{1,2} & \cdots & B_{1,m} & ~ e_1R \\
        B_{2,1} & R_2 & \cdots & B_{2,m} & ~ e_2R \\
        \vdots & \vdots & \ddots & \vdots & ~ \vdots \\
        B_{m,1} & B_{m,2} & \cdots & R_m & ~ e_mR \\
    \end{block}
\end{blockarray}
\]
This set with addition is an abelian group. To ensure it has the structure of a ring, one must impose axioms regarding products of coordinate bimodules, cf. [Han. 73a, Lemma 3.1]. A compact way used by many authors to ensure distributivity is to view multiplication of $B_{i,j}B_{j,k}$ as a $(R_i\text-R_k)$-bimodule \textit{product homomorphism} \[\varphi_{i,j,k}\colon B_{i,j}\otimes_{R_j} B_{j,k} \to B_{i,k}. \] It follows that multiplication by elements of $R_i$ coincides with its action on bimodules. The set of formal matrices then forms a ring as long as associativity holds for any triplet of bimodules different from $R_i$.

We can construct a formal matrix ring from any finite set of local rings and corresponding bimodules by setting all \textit{proper} product homomorphisms $\varphi_{i,j,k}$ with $i\neq j \neq k$ to zero. We call such rings \textit{trivial formal matrix rings}.

\subsection{Basic rings}\label{Secsempbasic}

Two rings are \textit{Morita equivalent} if their left (and hence also right) module categories are equivalent.  For a semiperfect ring $R$, there is (up to isomorphism) a unique canonical representative of the class of all rings Morita equivalent to $R$ called \textit{the basic ring for $R$}, as in the case of finite rings [McD. 74, Sec. X]. A basic ring for a QF ring is necessarily Frobenius.

\smallskip 

Recalling the ordering of a complete set of pairwise orthogonal idempotents $e_1, \dots, e_m$ from Sec. \ref{SecSempNot}, the basic ring for $R$ is isomorphic to  $End_R(e_1R\oplus \dots \oplus  e_nR)$. From a formal matrix ring, we can construct its basic ring as a leading principal submatrix of order $n$.

 Given a basic semiperfect ring of order $n$, all rings Morita equivalent to it are, up to isomorphism, determined by specifying the multiplicities $\mu_1, \mu_2, \dots, \mu_n$. From a basic formal matrix ring, we can then construct a Morita equivalent ring inductively by repeating rows and columns (block expansion). The multiplication maps on new coordinate bimodules are then induced by multiplication maps on their parent bimodules. 

 \[
\underbrace{
\begin{pmatrix}
  R_1 & B_{1,2} \\
  B_{2,1} & R_2
\end{pmatrix}
}_{\text{Basic Ring (Order } n=2)}
\xrightarrow{\text{Expansion } (\mu_1=2, \mu_2=1)}
\underbrace{
\left(
\begin{array}{cc|c}
  R_1 & R_1 & B_{1,2} \\
  R_1 & R_1 & B_{1,2} \\
  \hline
  B_{2,1} & B_{2,1} & R_2
\end{array}
\right)
}_{\text{Semiperfect Ring (Order } m=3)}
\]

Note that for $i, j\leq m$,  if $i,j$ are two indices from $I_k$, then $B_{i,j}$ is a copy of the local corner $R_k$. The multiplication by this coordinate bimodule is then induced by the corresponding action of $R_k$.
 
A property of a ring (or its category of modules) is called \textit{categorical} if it is preserved by a Morita equivalence.
 Many ring properties of a semiperfect ring $R$ studied in this article are categorical:

 semiperfect with a Nakayama permutation $\pi$, right self-injective, PF, QF, possessing a right (self-)duality [Xue 92, Cor. 4.6], right perfect, right artinian, finite, right Kasch, right QF-2, the inclusion $S_r\subseteq S_l$ and the essentiality of the right socle [NY 03, Lemma 3.19].

 On the other hand, combinatorial properties, like cardinality conditions, are not categorical. 
 
\subsection{Examples}\label{secsempex}

 In [Lam 99, Ex. 16.19(5)] T. Y. Lam writes: \say{\textit{Many books contain material on QF rings and Frobenius rings, but few of them offer any worked out examples of QF rings that are not Frobenius rings!}}, before presenting Nakayama's original example [Nak. 39, page 624]. The method of formal matrix ring representations provides an easy way to construct all QF rings that are Morita equivalent to a basic Frobenius ring. In particular, for a fixed basic Frobenius ring, we can construct all QF rings Morita equivalent to it. 

 \smallskip

Wood presents an example attributed to D. Benson [Wood 99, Ex. 1.4(1)]. Both this and Nakayama's example are subrings of rings of (classical) matrices over a field.  We worked out Wood's example in detail to illustrate the previous section. Let $R$ be a ring of matrices over a field $K$ of the form
\[
\left(
\begin{array}{cc|cc|cc}
     a&0 & b&0&0&0\\
     0&a & 0&b&c&0\\
     \hline
     d&0 & e&0&0&0\\
     0&d & 0&e&f&0\\
     \hline
     0&0 & 0&0&g&0\\
     h&0 & i&0&0&g
\end{array}
\right).
\]
 
 Then $diag(1,1,0,0,0,0)$, $diag(0,0,1,1,0,0)$, $diag(0,0,0,0,1,1)$ form a complete set of pairwise orthogonal primitive idempotents, and the corresponding local corners are isomorphic to $K$. We can then represent $R$ as $\left(
\begin{array}{cc|c}
  K&K&L\\
     K&K&L\\
  \hline
 L&L&K
\end{array}
\right)$, where $L$ is again a field $K$ viewed as a $(K\text-K)$-bimodule with multiplication defined as $LL=0$, whereas the multiplication of the form $KK$ is the standard multiplication of two elements from the field $K$.

 We reorder the idempotents (swapping indices 2 and 3) to maintain our notational preference of keeping the first $n$ idempotents ($n=2$ here) as a basic set of idempotents. We get the ring $\begin{pmatrix}
     K&L&K\\
     L&K&L\\
     K& L&K
 \end{pmatrix}$, with multiplicities $\mu_1=2$ and $\mu_2=1$, and $I_1=\{1, 3\}$, $I_2=\{2\}$, thus it is a QF ring that is not Frobenius.
 
 The basic ring associated with this ring is then $\begin{pmatrix}
     K&L\\
     L&K
 \end{pmatrix}$, which is a Frobenius ring whose Nakayama permutation is $(1~2)$. Note that it is isomorphic to [Wood 99, Ex. 1.4(3)].
 
 A QF ring whose Nakayama permutation is a cycle is Frobenius if and only if it is isomorphic to a (classical) matrix ring over its basic ring. In particular
\[
\left(
\begin{array}{cc|cc}
       K&L&K&L\\
     L&K&L&K\\
     \hline
     K&L&K&L\\
     L&K&L&K\\
\end{array}
\right)
\] with equal multiplicities $\mu_1=\mu_2=2$ is again Frobenius.

\subsection{Simple modules}\label{SecsempSimples}

 Green developed the theory of modules over formal matrix rings of order 2, generalising known results for generalised lower
triangular rings. We use his approach to present a canonical representation of simple modules over a semiperfect ring. This will provide a straightforward method for describing the socles of formal matrix rings.

\smallskip

Let $R$ be a semiperfect ring. Following [Gre. 82, Thm. 1.5], we can represent any right $R$-module $M$ as a row module $(M^1,\dots, M^m)$, where $M^i$ is a right $R_i$-module, together with multiplication maps $M^i\otimes_{R_i} B_{i,j}\to M^j$ defined in such a way that the ring action is associative. 

Let $N=(N^1,\dots, N^m)$ be a right $R$-module. Then $N$ is an $R$-submodule of $M$ if it is closed under ring actions, i.e., it holds that $N^iB_{i,j}\subseteq N^j$. In particular, $N^i$ is an $R_i$-submodule of $M^i$.

For $k\leq n$, we denote $U_k$ the unique simple right $R_k$-module and define a module $V_k$ such that its $i$-th coordinate is a copy of $U_k$ if $i\in I_k$ and zero otherwise. The multiplication is then a zero map, unless $B_{i,j}$ is a copy of $R_k$, i.e., if $i,j\in I_k$, in which case it is induced by the ring action of $R_k$.

\begin{Lemma}\label{LemmaSimplemodules}
    Modules $V_1, V_2, \dots V_n$ form, up to isomorphism, a complete set of simple $R$-modules.
\end{Lemma}
\begin{Proof}
    We prove that each of them is simple. Then, because they have different supports, they are pairwise non-isomorphic. In a semiperfect ring, simple modules are isomorphic to the tops of basic ideals; thus, there are $n$ of them, proving that our list is complete.

    Suppose that for some $k\leq n$, the module $V_k$ is not simple. Then there exists a nonzero proper submodule $W\subseteq V_k$. Because all nonzero coordinates are simple modules, we get that there are two indices $i,j\in I_k$ such that $W^i=U_k$ and $W^j=0$. But then $W^iB_{i,j}=U_kR_k=U_k\neq 0$, which is a contradiction with $W$ being a submodule.
    \end{Proof}
The left simple modules are described dually.

\begin{Ex}
    Compare this with the description of simple modules for a (classical) matrix ring $M_m(K)$ of order $m$ over a field $K$. We get that $n=1$, $\mu_1=m$ and $I_1=\{1, 2, \dots, m\}$. Then $U_1=K$ and we have one unique simple right module, namely $(K, K, \dots, K)$.
\end{Ex}

\subsection{Morita Duality}\label{SecDual}

In this section, we survey the basics of Morita duality following [Xue 92], which serves as the primary reference for the results presented here. Although we do not assume prior knowledge of Morita duality in the sequel, the close connection between our results and duality theory makes a brief overview beneficial. In particular, the concepts of \textit{i-pairs} and \textit{Nakayama pairs} provide insight into the structure of rings with Nakayama permutations."

\smallskip

Let $R$ and $S$ be two rings and $_SU_R$ an $(S\text-R)$-bimodule. Then we can define a pair of contravariant functors \[
\xymatrix{
    S\text{-Mod} \ar@<0.5ex>[rr]^{Hom_S(-, _SU)} 
    & & Mod\text{-R} \ar@<0.5ex>[ll]^{Hom_R(-, U_R)}
}\]

Viewing $R$ as an $(R\text-R)$-bimodule gives a functor $Hom_R(-, R)$, often called  \textit{duality} in the literature, even in cases where it is not a Morita duality (e.g., [Die. 58]).

A \textit{Morita Duality} is a contravariant equivalence between two subcategories of  $S\text-Mod$ and $Mod\text-R$ that are closed under submodules and factor modules and contain $_SS$ and $R_R$, respectively. Such a duality, if it exists, is always given by an $(S\text-R)$-bimodule $_SU_R$, which is a faithful injective cogenerator. In particular, $U$ contains an isomorphic copy of the injective envelopes of all isomorphism types of simple right $R$-modules (left $S$-modules).  We will refer to such a bimodule as a \textit{duality context}. 

Over artinian rings, the duality context coincides with \textit{quasi-Frobenius bimodules} introduced by G. Azumaya [Azu. 59]. This name is still in use in algebraic coding theory, where finite quasi-Frobenius bimodules are used as an alphabet for codes. See [GNW 04] for some recent results and an overview of the literature.

If  $R$ has a Morita duality, then both $R_R$ and the minimal injective cogenerator in $Mod\text-R$ are \textit{linearly compact}. The induced duality is precisely that of the category of right linearly compact $R$-modules, independent of the specific duality context. For more on linear compactness, see [Xue 92, Ch. 1.3]. Note that any artinian module is linearly compact. In particular, over a right perfect ring, linearly compact left modules are precisely artinian modules [Xue 92, Prop. 18.3]. If the ring is perfect, artinian modules have a finite length [Xue 92, Prop. 18.2]. Thus, over artinian rings, linearly compact modules are precisely modules of finite length. Over finite rings, these are the finite modules.

\begin{Ex}
 Artin algebras have a self-duality. A ring $R$ is called an \textit{artin algebra} if it is artinian and finitely generated over its centre. Both fin. dim. algebras over a field and finite rings are artin algebras. 
\end{Ex}

Artinian rings with duality, which is not a self-duality, exist [Xue 92, Ex. 12.7], as do artinian rings without duality [Xue 92, Remark 2.9].

\smallskip 

Note that if two rings $S$ and $S'$ are Morita duals of $R$, then $S$ and $S'$ are Morita equivalent, as the composition of dualities induces an equivalence between the categories of their linearly compact left modules. The equivalence of categories of linearly compact modules implies Morita equivalence, and the proof mirrors the argument showing that the equivalence of f.g. modules implies Morita equivalence [AF 92, Exc. 22.4].

\begin{Prop}\label{PropDual}
    Let $S$ and $S'$ be two Morita duals of $R$. 

    Then $S$ and $S'$ are Morita equivalent.
\end{Prop}
\begin{Proof}
    Let us denote $S\text-lc$ and $S'\text-lc$ the respective categories of left $S$ and $S'$-modules. Composing their respective dualities with $lc\text-R$ we obtain an equivalence of categories \[
\xymatrix{
   S\text-lc \ar@<0.5ex>[rr]^{F} 
    & & S'\text-lc. \ar@<0.5ex>[ll]^{G}
}\]
    
    First, observe that if $P:=F(_SS)$ is a f.g. projective generator, the result follows from [Lam 99, Thm. 18.24], because we have ring isomorphisms 
    \[S\cong End(_SS)\cong End(F(S)).\]
    To see the second isomorphism, observe that an additive category equivalence commutes with $Hom(-, -)$ so $End(_SS)\cong End(F(S))$ as groups. Since $F$ is a functor, it preserves composition of morphisms, so it respects multiplication, showing that this is indeed a ring isomorphism.

    Both $F(S)$ and $G(S')$ are f.g. because $_SS$ and $_{S'}S'$ are f.g., and it is a categorical property [Lam 99, 18.2]. Because linearly compact modules are closed under submodules, we can modify the classical argument that shows that being f.g. is a categorical property to show that being linearly compact is a categorical property too.
    
    Since $S$ is projective in $S\text-Mod$, it is a projective object in $S\text-lc$, so $F(S)$ is projective in $S'\text-lc$, because category equivalence preserves projective objects. Similarly, $G(S')$ is projective in $S\text-lc$.
    
    Because $G(S')$ is finitely generated, there is $n$ such that there exists an epimorphism \[\pi\colon S^n\twoheadrightarrow G(S').\] As $S$ is linearly compact, so is $S^n$, so this $\pi$ is an epimorphism in $S\text-lc$ where $G(S')$ is projective, so $\pi$ splits. 

    An equivalence preserves split epimorphisms, so  $F(\pi)$ is a split epimorphism 
    \[F(S)^n\cong F(S^n) \overset{F(\pi)}{\twoheadrightarrow} FG(S')\cong S', \]
    so $S'$ is a direct summand of a finite direct power of $F(S)$, so $F(S)$ is a generator in $S'\text-Mod$ [Lam 99, Thm. 18.8(3)].

    Because $F(S)$ is finitely generated, there is an epimorphism \[p\colon (S')^m\twoheadrightarrow F(S).\] Once again, $p$ is an epimorphism in $S'\text-lc$ so it splits because $F(S)$ is a projective object in this category. However, this means that $F(S)$ is a direct summand of a free $S'$-module; hence, it is a projective module in $ S'\text {- Mod} $.    
\end{Proof}

Rings with Morita duality are necessarily semiperfect [Xue 92, Thm. 2.7]. Recalling the notation of coordinate bimodules $B_{i,j}$ (Sec. \ref{Secsempmatrix}), we have the following characterisation of rings with Morita duality in terms of their formal matrix representation: 
\begin{Prop}[Xue 92, Thm. 4.12]\label{PropDuality}
Let $R$ be a semiperfect ring.

Then $R$ has a right Morita duality if and only if each local corner $R_i$ has a right Morita duality given by an $R_i$-module $U_i$, and for all $i, j\leq m$, the bimodules $B_{i,j}$ and $Hom_{R_i}(B_{j,i}, U_i)$ are linearly compact as right $R_j$-modules.    
\end{Prop}

QF rings are the class of artinian rings when the standard duality $Hom_R(-, R)$ is a Morita self-duality. A ring such that $Hom_R(-,R)$ is a Morita self-duality is called a \textit{cogenerator ring}. Cogenerator rings coincide with PF rings.

\begin{Prop}[Kato 67, Thms. 7.,  10.]\label{PropKato}
 Let $R$ be a ring. TFAE

 (1) $R$ is a PF ring.

 (2) $R$ is a self-injective Kasch ring.

 (3) $R$ is a self-injective ring with a Nakayama permutation.

 (4) $R$ is a self-injective D-ring. 
\end{Prop}
A ring is called D-ring if it satisfies the Double annihilator property (Def. \ref{DefD}).

Another prominent class of duality contexts consists of $(R_i\text-R_{\pi(i)})$-bimodules of the form $e_iRe_{\pi(i)}$ in PF rings, strengthening Prop. \ref{PropDuality}. The artinian case was proved by K. Fuller [Ful. 69, Thm. 3.1] and generalised by B. Roux [Roux 71, Sec. 3]. This observation was used to give a description of formal matrix representations of PF rings by Müller (Sec. \ref{SecMueller}).

These results are based on the concept of an \textit{i-pair}, which generalises some properties of Nakayama permutations. The term $i$-pair originates from Fuller [Ful. 69, Thm. 3.1], although his definition of a \textit{pair} differs. For an overview of the results and literature on i-pairs, see [BO 09]. 

\begin{Def}[i-pairs]
Let $R$ be a semiperfect ring and $e, f\in R$ primitive idempotents.

The couple  $(eR, Rf)$ is an \emph{$i$-pair} if and only if $soc(eR)\cong top(fR)$ and $soc(Rf)\cong top(Re)$.

If both $eR$ and $Rf$ have essential socles, $(eR, Rf)$ is called a \emph{Nakayama pair}.
\end{Def}
 In a basic ring with a Nakayama permutation $\pi$, all $i$-pairs are of the form $(e_iR, Re_{\pi(i)})$. Xue coined the term Nakayama pair. He realised that Fuller's results apply to an arbitrary semiperfect ring as long as $eR$ and $fR$ have essential socles.

\begin{Prop}[Xue 97, Cor. 11, Thm. 7]\label{PropIpairsperfect}
Let $R$ be a ring with a Nakayama permutation $\pi$ and let $k\leq n$ be such that the basic ideals $e_kR$ and $Re_{\pi(k)}$ have essential socles. TFAE

(1) Both $e_kR\in e_kRe_k\text-mod$ and $Re_{\pi(k)}\in mod\text- e_{\pi(k)}Re_{\pi(k)}$ are linearly compact.

(2) $e_kR_R$ is injective, and $Re_{\pi(k)}\in mod\text- e_{\pi(k)}Re_{\pi(k)}$ is linearly compact.

(3) Both $e_kR_R$ and $_R Re_{\pi(k)}$ are injective.

(4) $_R Re_{\pi(k)}$ is injective, and $e_kR\in e_kRe_k\text-mod$ is linearly compact.

Furthermore, if $R$ is perfect, it is enough to assume that only one of the modules in (1) is linearly compact.     
\end{Prop}
It follows that the self-injectivity in Prop. \ref{PropKato}(3) can be replaced by linear compactness or an analogous condition using Prop. \ref{PropIpairsperfect} parts (2) or (4) [Xue 97, Cor. 9].

Combining Props. \ref{PropDuality} and \ref{PropIpairsperfect}, we get the following possibly new corollary:
\begin{Cor}\label{CorDualPF}
Let $R$ be a ring with a Nakayama permutation and essential socles.

Then $R$ has a Morita self-duality if and only if it is a PF ring.   

In particular, a D-ring has Morita self-duality if and only if it is a PF ring.
\end{Cor}
\begin{Proof}
PF rings are known to have Morita self-duality. For the converse, assume that $R$ has a Morita self-duality. Then $R$ is linearly compact, and by Prop. \ref{PropDuality}, all bimodules $e_iRe_j$ are linearly compact. A right (left) $R$-module $M$ is linearly compact if and only if right (left) $e_iRe_i$-modules $Me_i$ ($e_iM$) are linearly compact for all $i\leq m$. Thus, all primitive ideals are linearly compact as $R$-modules. 

    By assumption, all primitive ideals have essential socles, so $R$ satisfies  Prop. \ref{PropIpairsperfect}(1), and hence, by Prop. \ref{PropIpairsperfect}(3), $R$ is a self-injective ring.  

    The \say{in particular} part follows because D-rings have a Nakayama permutation and essential socles (Thm. \ref{ThmD}).
    \end{Proof}
    By the \say{furthermore} part of Prop. \ref{PropIpairsperfect}, if $R$ is perfect, it is enough to assume any Morita duality, not necessarily self-duality.
\begin{remark}
    It is known that D-rings need not be self-injective [HN 85,  Exs. 6.1, 6.2]. So, we see that lattice self-duality does not imply Morita self-duality. Because rings with Morita duality are described purely by structural properties (Prop. \ref{PropDuality}), it follows that rings with self-Morita duality that are not D-rings exist. 
\end{remark}
Corollary \ref{CorDualPF} characterises the class of PF rings as the intersection of two classes of rings defined by self-duality: lattice self-duality and Morita self-duality. We now show that for simple modules, the two types of duals coincide. 
\begin{Prop}
Let $R$ be a  PF ring and $\pi$ be its Nakayama permutation.

Then the dual of a simple right ideal $T_k=soc(e_kR)$ is isomorphic to $R/l(T_k)\cong V'_{\pi^{-1}(k)}$.
\end{Prop}
By $l(T_k)$ we denote the left annihilator of $T_k$. As PF rings are Kasch, this proposition describes the duals of all simple modules. The case for QF rings is proved in [Wood 99, Cor. 2.5].
\begin{Proof}
Because $R$ is PF, $_RR_R$ is the corresponding self-duality context, and hence the dual of $T_k$ is equal to $Hom_R(T_k, R)$. Because $R$ is right self-injective,  elements of  $Hom_R(T_k, R)$ can be represented as left multiplications by elements of $R$. Formally, we have an epimorphism of left modules that maps $r\in R$ to a module homomorphism $r\cdot \colon T_k\to R$. The kernel of this map is precisely the elements that vanish on $T_k$. Hence, by the isomorphism theorem, $R/l(T_k)\cong Hom_R(T_k, R)$ as left modules.

The isomorphism $R/l(T_k)\cong V'_{\pi^{-1}(k)}$ is a classical result due to Nakayama, which we recall in Lemma \ref{LemmaAnnofsimple}.
    \end{Proof}

\section{Definition of a Nakayama permutation}\label{SecDef}

This section defines the Nakayama permutation for a general semiperfect ring. A ring with a Nakayama permutation is necessarily Kasch, and all primitive ideals have simple socles. We show that as long as $S_r$ or $S_l$ is essential, this observation actually characterises rings with Nakayama permutation (Thm. \ref{PropQFKaschEss}). If such a ring is Noetherian, it is necessarily a QF ring (Thm. \ref{ThmNY}). We then recall Müller's description of PF rings (Sec. \ref{SecMueller}) and construct several counterexamples showing that our results and Müller's cannot be further generalised (Sec. \ref{SecDefcounter}).

\smallskip 

The Nakayama permutation was originally formulated for fin. dim. algebras over fields in Nakayama's seminal paper [Nak. 39, Lemma 2]. Generalising his results to the artinian setting in the subsequent paper [Nak. 41a, Chap. 2], he used the Nakayama permutation as the definition of a QF ring.

For any $n>1$, we denote the set of permutations of $\{1, 2, \dots, n\}$ by $S_n$.
\begin{Def}[Nakayama permutation]\label{DefNak}
    Let $R$ be a semiperfect ring, $e_1, \dots, e_n$ a basic set of idempotents, and let $\pi\in S_n$. 

    Then $\pi$ is called a \emph{Nakayama permutation of $R$} if, for each $k\leq n$
    \begin{gather}
        soc(e_kR)\cong top(e_{\pi(k)}R)\\
        soc(Re_{\pi(k)})\cong top(Re_{k}).
    \end{gather}
\end{Def}
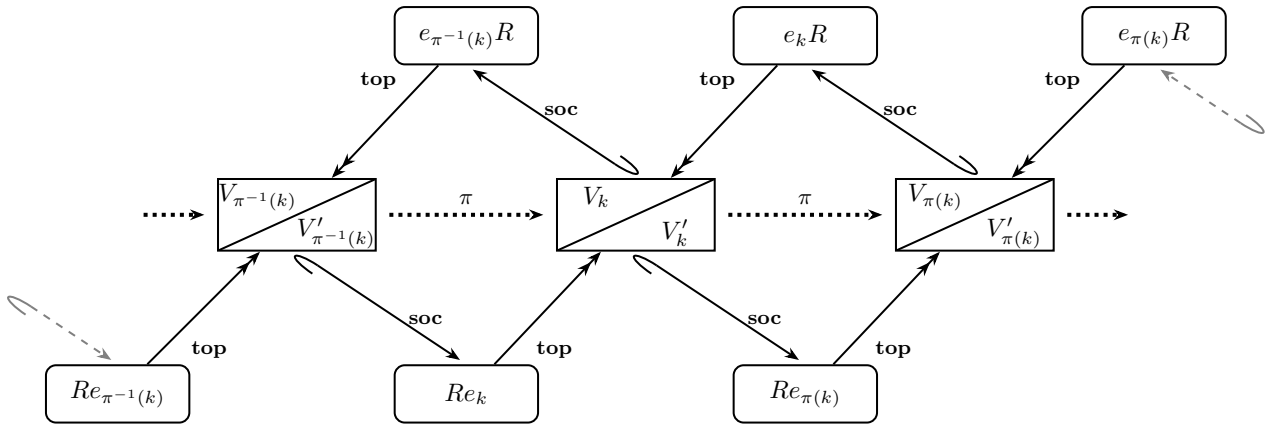
\begin{figure}[ht]
    \centering
    \resizebox{\textwidth}{!}{
    \begin{tikzpicture}[
        thick,
        >={Stealth[length=2mm]}, 
        node distance=1.5cm and 2.5cm,
        splitblock/.style n args={2}{
            draw=black,
            rectangle,
            minimum width=2.2cm,
            minimum height=1.0cm,
            fill=white,
            font=\small,
            path picture={
                \draw[black] (path picture bounding box.south west) -- (path picture bounding box.north east);
                \node[anchor=center] at ($(path picture bounding box.west)!0.5!(path picture bounding box.north)$) {#1};
                \node[anchor=center] at ($(path picture bounding box.east)!0.5!(path picture bounding box.south)$) {#2};
            }
        },
        ideal/.style={
            draw=black,
            rectangle,
            rounded corners,
            minimum width=2.0cm,
            minimum height=0.8cm,
            align=center,
            fill=white,
            font=\small
        },
        progression/.style={
            ->,
            dotted,
            ultra thick,
            shorten >=5pt,
            shorten <=5pt,
            font=\footnotesize
        },
        socmap/.style={
    {Hooks[left, length=3.5mm, width=3mm]}->,  
    shorten >=3pt,
    thick,
    font=\scriptsize\bfseries
},
        topmap/.style={
            ->>,
            font=\scriptsize\bfseries
        }
    ]

    \node[splitblock={$V_k$}{$V'_k$}] (Vk) at (0,0) {};
    \node[splitblock={$V_{\pi^{-1}(k)}$}{$V'_{\pi^{-1}(k)}$}, left=of Vk] (Vprev) {};
    \node[splitblock={$V_{\pi(k)}$}{$V'_{\pi(k)}$}, right=of Vk] (Vnext) {};

    \node[ideal] (ePrevR) at ($(Vprev)!0.5!(Vk) + (0, 2.5)$) {$e_{\pi^{-1}(k)}R$};
    \node[ideal] (ekR) at ($(Vk)!0.5!(Vnext) + (0, 2.5)$) {$e_kR$};
    \node[ideal] (eNextR) at ($(Vnext) + (2.5, 2.5)$) {$e_{\pi(k)}R$};

    \node[ideal] (RePrev) at ($(Vprev) + (-2.5, -2.5)$) {$Re_{\pi^{-1}(k)}$};
    \node[ideal] (Rek) at ($(Vprev)!0.5!(Vk) + (0, -2.5)$) {$Re_k$};
    \node[ideal] (ReNext) at ($(Vk)!0.5!(Vnext) + (0, -2.5)$) {$Re_{\pi(k)}$};

    
    \draw[progression] ($(Vprev.west)+(-1.2,0)$) -- (Vprev.west);
    \draw[progression] (Vprev.east) -- node[above] {$\pi$} (Vk.west);
    \draw[progression] (Vk.east) -- node[above] {$\pi$} (Vnext.west);
    \draw[progression] (Vnext.east) -- ($(Vnext.east)+(1.2,0)$);

    

    \draw[topmap] (ePrevR) -- node[above left, pos=0.3] {top} (Vprev);
    \draw[topmap] (RePrev) -- node[below right, pos=0.3] {top} (Vprev);

    \draw[topmap] (ekR) -- node[above left, pos=0.3] {top} (Vk);
    \draw[topmap] (Rek) -- node[below right, pos=0.3] {top} (Vk);

    \draw[topmap] (eNextR) -- node[above left, pos=0.3] {top} (Vnext);
    \draw[topmap] (ReNext) -- node[below right, pos=0.3] {top} (Vnext);

    
    \draw[socmap] (Vk.north) -- node[right, pos=0.6] {soc} (ePrevR.south);
    \draw[socmap] (Vnext.north) -- node[right, pos=0.6] {soc} (ekR.south);
    \draw[socmap, dashed, gray] ($(eNextR.south)+(1.5, -1.0)$) -- node[right, pos=0.6] {} (eNextR.south);

    \draw[socmap, dashed, gray] ($(RePrev.north)+(-1.5, 1.0)$) -- node[right, pos=0.6] {} (RePrev.north);
    \draw[socmap] (Vprev.south) -- node[right, pos=0.6] {soc} (Rek.north);
    \draw[socmap] (Vk.south) -- node[right, pos=0.6] {soc} (ReNext.north);

    \end{tikzpicture}
    }
    \caption{Visualizing the Nakayama permutation $\pi$. }
    \label{FigNakayamaStraightDiag}
\end{figure}
Recall that all primitive right (left) ideals are isomorphic to some basic right (left) ideal $e_kR$ ($Re_k$), showing $top(e_{\pi(k)}R)\cong V_{\pi(k)}$. Thus, all primitive right  (left) ideals have simple socles. J. Wood calls this \say{\textit{the real content hidden in [the definition of Nakayama permutation]}} in [Wood 99, page 557].  R. M. Thrall, who pioneered the study of algebras with this property, coined the term \textit{QF-2 algebra} [Thr. 48]. 
The name QF-2 refers to it being one of the three generalisations of QF-rings considered by Thrall.

\begin{Def}[QF-2 rings]
    Let $R$ be a semiperfect ring. 

    Then $R$ is called a \emph{right (left) QF-2} ring if the socle of every right (left) primitive ideal is simple. 
\end{Def}
 Recalling the description of projective modules over semiperfect rings (Sec. \ref{SecSempNot}), $R$ is right QF-2 if and only if all indecomposable f.g. projective right modules have simple socles.

 Furthermore, Definition \ref{DefNak} also implies the following property, which we discuss in detail in Sec. \ref{SecDKasch}: 
\begin{Def}[Kasch rings]
    Let $R$ be a semiperfect ring.

    Then $R$ is \emph{right (left) Kasch} if every simple right (left) module is isomorphic to a right (left) ideal of $R$.
\end{Def}
\smallskip

 In the remark adjoining the definition of QF rings, Nakayama observes that artinian left QF-2 rings satisfying right socle condition (1) possess a Nakayama permutation. F. Kasch in [Kas. 82, Thm. 13.4.2] proves that for artinian rings, the left and right versions of the Nakayama permutation imply the Nakayama permutation. We note that this proof applies verbatim to the more general setting of semiperfect rings with essential socles. He also proves that an artinian ring has a Nakayama permutation if and only if it is a QF-2 ring and its socles coincide.
 
 In general, rings with a Nakayama permutation do not necessarily have coinciding socles (Ex.  \ref{ExSocsDifer1}), and Kasch QF-2 rings need not have a Nakayama permutation  (Exs. \ref{ExQFKasch},  \ref{ExQFKasch2}). However, the structure of such rings is not well understood.  On the other hand, QF-2 rings with coinciding socles have a Nakayama permutation (Prop.  \ref{PropSocsCoincide}).

\subsection{Kasch QF-2 rings}\label{SecDefKasch}

This section begins by isolating and examining the consequences of a single Nakayama permutation condition (Def. \ref{DefNak}(1)), demonstrating its equivalence with the right Kasch property for semiperfect right QF-2 rings (Lemma \ref{Lemma(1)}).  We then address a critical distinction: while both the left and right versions of the permutation may exist, they are not guaranteed to be inverses of each other. We show that it is enough for the right and left socles to coincide (Prop. \ref{PropSocsCoincide}), a condition that is guaranteed if one of them is essential  (Thm. \ref{PropQFKaschEss}). If left and right socles coincide in a ring with a Nakayama permutation, then simple modules linked by a Nakayama permutation have isomorphic rings of endomorphisms (Prop. \ref{PropIsomfields}).

\smallskip 

When discussing only one condition in isolation, say (1), there are several equivalent module-theoretic reformulations.
\begin{Lemma}\label{Lemma(1)}
    Let $R$ be a semiperfect right QF-2 ring. The following are equivalent 
    
     (1) There exists a permutation $\pi\in S_n$ that satisfies Def. \ref{DefNak}(1),

    (2) $R$ is right Kasch.

    (3) Primitive right ideals are, up to isomorphism, determined by their socles. 

    (4) Projective right modules are, up to isomorphism, determined by their socles. 
\end{Lemma}
\begin{Proof}
    (1)$\implies$(2) Basic ideals are pairwise nonisomorphic projective modules. Thus, their tops are pairwise nonisomorphic. Since they all have simple socles that are paired to these tops by condition (1), there are at least $n$ nonisomorphic simple right modules contained in $R$; hence, all of them.

    (2)$\implies$(3) Suppose two non-isomorphic primitive ideals had isomorphic socles. Since there are $n$ isomorphism classes of primitive ideals and $n$ classes of simple modules, the Pigeonhole Principle implies that the map $e_i R \mapsto soc(e_i R)$ would miss an isomorphism type of simple module. Thus, $S_r$ would not contain that missing type, as any simple right ideal is isomorphic to a socle of a right primitive ideal (Lemma \ref{LemmaSoceR}).

    (3)$\implies$(4) By [Mül. 70, Thm. 3], right projective modules of a semiperfect ring are direct sums of isomorphic copies of primitive right ideals. Thus, they can be, up to isomorphism, characterised by the number of copies of each basic ideal. Using (3), they are thus determined by the number of copies of each (isomorphism type of) simple right module in their socle. 

    (4)$\implies$(1) Applying (4) to basic ideals, each of them has a pairwise nonisomorphic simple socle. 
\end{Proof}

However, a QF-2 Kasch ring is not necessarily a ring with a Nakayama permutation. The left and right permutations, even if they exist, may not be inverses of each other (Exs. \ref{ExQFKasch}, \ref{ExQFKasch2}).

From the perspective of a basic formal matrix ring, the necessity of coinciding socles is evident. The general proof is provided in the next proposition. Coinciding socles turn out to be a rather natural condition for such rings (Thms. \ref{PropQFKaschEss}, \ref{ThmAnn1}).

The idea for the following proof is based on [HN 85, Thm. 4.2], though the result may predate this reference.
\begin{Prop}\label{PropSocsCoincide}
    Let $R$ be a semiperfect QF-2 ring with coinciding socles.

    Then $R$ has a Nakayama permutation.
\end{Prop}
\begin{Proof}
Take $k\leq n$, then $e_ksoc(R)=soc(e_kR)$ is a simple module because the ring is right QF-2. By the right version of [CR 64, Thm. 54.12] (cf. [McD. 74, Thm. XI.1]) we have that $e_ksoc(R)\cong top(e_lR)$ if and only if $e_ksoc(R)e_l\neq 0$. Because $e_1, \dots, e_n$ is a basic set of idempotents, there is a unique $l\leq n$ such that the simple right module $e_ksoc(R)$ is isomorphic to $top(e_lR)$. We denote this by $\pi(k):=l$. Then $\pi$ is a permutation of $\{1, 2, \dots, n\}$ and by the dual version of the above argument, we show that it is actually a Nakayama permutation.    
\end{Proof}

We now aim to show that a Kasch QF-2 ring with an essential right socle has coinciding socles and hence a Nakayama permutation. 

\begin{Lemma}
    Let $R$ be a basic semiperfect ring of order $n$. Further assume it is left QF-2, left Kasch, and its right socle is square-free essential right ideal.

    Then the left and right socles of $R$ coincide.
\end{Lemma}
\begin{Proof}
First, observe that $S_l$ is square-free. It is generated by simple modules of the form  $soc(e_kR)$. These are pairwise nonisomorphic by the Kasch condition.  By Lemma \ref{LemmaUrelated}, this means that $R$ has precisely $n$ simple left ideals. In particular, since every homogeneous component of $S_l$ is a two-sided ideal and is simple as a left ideal, every simple left ideal is also a right ideal.

   Let $T'$ be an ideal that is simple as a left ideal. As noted above, $T'$ is also a right ideal.  Because $S_r$ is an essential right ideal, the intersection $S_r\cap T'$ is nonzero. On the other hand, both $S_r$ and $T'$ are left ideals; their intersection is also a left ideal. But as a left ideal, $T'$ is simple, so $T'\subseteq S_r$, proving $S_l\subseteq S_r$.

    Because simple left ideals are contained in the right socle, they are semisimple as right modules. But because $S_r$ is square-free, we once again infer that there are exactly $n$ simple right ideals. Thus, all simple left ideals are simple as right ideals and constitute the whole $S_r$.
\end{Proof}
\begin{Thm}\label{PropQFKaschEss}
    Let $R$ be a QF-2 Kasch ring with an essential right socle.  

    Then $S_l=S_r$ and $R$ has a Nakayama permutation.

    In particular, any semiperfect ring with a Nakayama permutation and essential right socle has coinciding socles. 
\end{Thm}
\begin{Proof}
    Since $R$ is QF-2, by Proposition \ref{PropSocsCoincide}, it is enough to prove that its socles coincide. Now, this and the assumptions of the theorem are Morita-invariant properties, so we can assume that the ring is basic. By the previous Lemma, it is enough to observe that its right socle is square-free. However, it follows from being right QF-2, right Kasch, and from Lemma \ref{LemmaUrelated}. 
\end{Proof}

The observation that a ring with a Nakayama permutation and an essential right socle must have coinciding socles is explicitly shown in the proof of [NY 03, Thm. 3.21]  and follows from [Nor. 75, Prop. 2.2.4, Thm. 2.2.9].

We conclude this section with the following result, showing that simple modules linked by the Nakayama permutation have isomorphic rings of endomorphisms.  The special case for PF rings also follows from Müller's description of PF rings (Ex. \ref{ExMueler}), but was not explicitly formulated.
\begin{Prop}[Kra. 25, Prop. 18]\label{PropIsomfields}
    Let $R$ be a semiperfect ring with a Nakayama permutation $\pi$ and coinciding socles, and let $k\leq n$.

    Then $\m_k\cong \m_{\pi(k)}$ as skew fields.
\end{Prop}
This isomorphism is fundamental to the structure of rings with a Nakayama permutation. It ensures that the homogeneous components of the top and socle corresponding to $k$ and $\pi(k)$ share the same \say{geometry} (skew fields). As we will see in Sec. \ref{SecDouble}, this geometric invariance allows the lattice anti-isomorphism of the socle and the top to survive Morita equivalence, even when differing multiplicities prevent the top and socle from being isomorphic as modules

\subsection{Noetherian case}\label{SecDefNoether}

It is a recurring theme in the study of Nakayama permutations that right Noetherian conditions often imply that the ring is right artinian.  As a prominent example, a left or right Noetherian right self-injective ring is a QF ring [Lam 99, Thm. 15.1]. 

In view of Theorem \ref{PropQFKaschEss}, for rings with a Nakayama permutation, many results on artinian rings can be generalised using the following result of J. Chen, N. Ding, and M. Yousif.
\begin{Thm}[CDY 04, Thm. 2.5]
Let $R$ be a right Noetherian ring with an essential right socle such that $S_r\subseteq S_l$.

Then $R$ is right artinian.    
\end{Thm}
In particular, a right Noetherian ring with a Nakayama permutation and an essential right socle satisfies these conditions (Thm. \ref{PropQFKaschEss}) and is therefore right artinian. Consequently, such a ring is QF if and only if it is a minannihilator ring, which generalises a classical result by H. H. Storrer (Prop. \ref{PropStorrer}).

\subsection{Müller's description of PF rings} \label{SecMueller}

Rings with cyclic Nakayama permutations are the building blocks for all rings possessing a Nakayama permutation. For QF rings, this was observed by T. A. Hannula [Han. 73a,  Prop. 2.1], who also described how to glue QF rings together in terms of Morita contexts [Han. 73a, Prop. 3.3, Thm. 3.6]. The case for PF rings was proved by Müller in [Mül. 74, Lemma 5], where he also described the matrix representation of PF rings with a cyclic Nakayama permutation (Ex. \ref{ExMueler}). Recently, Müller's results have been generalised to rings whose right socle is essential [Kra. 25, Sec. 4]. A combinatorial method of glueing rings with cyclic Nakayama permutation together is described in [Kra. 25, Sec. 5].

\begin{Thm}[Kra. 25, Thm. 17]\label{ThmCycle}
      Let $R$ be a basic semiperfect ring with a Nakayama permutation $\pi$ and essential socles. Furthermore, let $I\subseteq \{1,2,\dots, n\}$ be a nonempty set of indices such that if $i\in I$ then $\pi(i)\in I$ and set $e:=\sum_{i\in I}e_i$. 

    Then $eRe$ is a ring with a Nakayama permutation $\pi_{\restriction_I}$ and essential socles.

    Furthermore, if $R$ is PF/QF, then so is $eRe$.
\end{Thm}

 In the remainder of this section, let $R$ be a PF ring, whose Nakayama permutation is an $n$-cycle. Following [Mül. 74, Thm. 4],  rings $R_1, \dots, R_n$ can be local corners of a ring whose Nakayama permutation is $(1~2~\dots~n)$ if and only if they form a \textit{cycle of dualities}: that is, for each $k\leq n$, there exists a duality between left (linearly compact) $R_{k+1}$-modules and right (linearly compact) $R_k$-modules (with indices taken modulo $n$).  It seems to be an open question whether a cycle of dualities can consist of non-isomorphic rings, even though local rings (even artinian) with a duality that does not possess self-duality exist [Xue 89].

Following Prop. \ref{PropDual},  if any one of the local corners $R_k$ has a self-duality, then all local corners are isomorphic. This, in particular, includes the case when $R$ is finite. In that case, $R$ has the following representation.
 \begin{Ex}[Mül 74, Thm.4]\label{ExMueler}
     Let $R$ be a basic finite Frobenius ring whose Nakayama permutation is an $n$-cycle.
     
     Then there exists a finite local ring $S$ such that all local corners of $R$ are isomorphic to $S$. Let $E$ denote an injective envelope of the unique simple $S$-module.
     
     Then there exists a collection of finite $S\text-S$-bimodules $B_1, \dots, B_h$, where $h:=\lfloor n/2\rfloor -1$. Futhermore, if $n$ is odd, let $C$ be an $S\text-S$-bimodule such that $C=Hom_S(C, S)$. Let $B_i^*:=Hom_S(B_i, S)$ for any $i\leq h$. 

     Then we can represent $R$ as in the following figure. \end{Ex}
\begin{figure}[ht]
    \centering
    \begin{tikzpicture}[
        x=0.6cm, y=0.6cm,
        cell/.style={draw, rectangle, minimum size=0.6cm, align=center, thick, font=\tiny, anchor=center},
        gap/.style={draw=none, rectangle, minimum size=0.6cm, align=center, font=\small, anchor=center},
        S_style/.style={fill=gray!10, font=\bfseries\tiny},
        E_style/.style={fill=gray!30},
        B_style/.style={fill=white},
        C_style/.style={fill=black!70, text=white, font=\bfseries\tiny}
    ]

    \foreach \myRow in {0,...,10} {
        \foreach \myCol in {0,...,10} {
            
            \pgfmathparse{int(mod(\myCol - \myRow + 11, 11))}
            \let\myDiag\pgfmathresult

            \ifcase\myDiag
                \node[cell, S_style] at (\myCol, -\myRow) {$S$};     
            \or
                \node[cell, E_style] at (\myCol, -\myRow) {$E$};     
            \or
                \node[cell, B_style] at (\myCol, -\myRow) {$B_1$};   
            \or
                \node[gap] at (\myCol, -\myRow) {$\cdot$};           
            \or
                \node[gap] at (\myCol, -\myRow) {$\cdot$};           
            \or
                \node[cell, B_style] at (\myCol, -\myRow) {$B_h$};   
            \or
                \node[cell, C_style] at (\myCol, -\myRow) {$C$};     
            \or
                \node[cell, B_style] at (\myCol, -\myRow) {$B_h^*$}; 
            \or
                \node[gap] at (\myCol, -\myRow) {$\cdot$};           
            \or
                \node[gap] at (\myCol, -\myRow) {$\cdot$};           
            \or
                \node[cell, B_style] at (\myCol, -\myRow) {$B_1^*$}; 
            \fi
        }
    }
    \end{tikzpicture}
    \label{FigMuel}
        \caption{Representation of a basic finite Frobenius ring whose Nakayama permutation is a cycle.}
\end{figure}
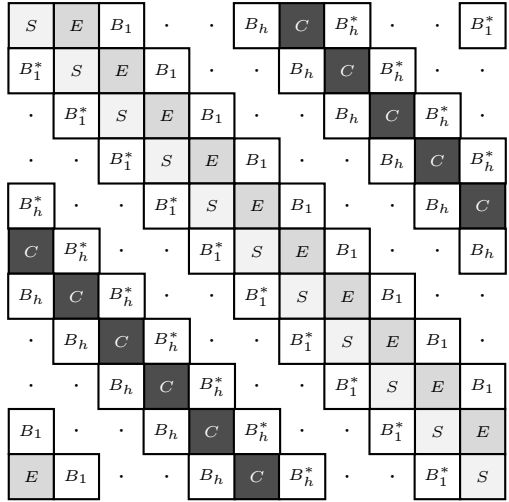
We remark that it is possible for bimodules in the same diagonal to not necessarily be the same, but only isomorphic, and the isomorphism problem for these representations remains open. Similarly, the problem of characterising the multiplication homomorphisms is unsolved, and only partial results are known [Han. 73b, Thm.11; Mül. 74, Lemma 3(3); Kra. 25, Thm. 13(b), (b')].

\subsection{Counterexamples}\label{SecDefcounter}

The following examples show that the results of this section rely crucially on the essentiality of the right socle. In particular, Kasch QF-2 rings need not have a Nakayama permutation (Exs. \ref{ExQFKasch}, \ref{ExQFKasch2}). A ring with a Nakayama permutation need not have coinciding socles if neither of the socles is essential (Ex. \ref{ExSocsDifer1}, and the symmetries observed by Müller (ex. \ref{ExMueler}), cf. [Kra. 25, Prop. 25]) need not hold (Ex. \ref{ExnotEss}).

\smallskip 

Recall that the top of a local ring is a \textit{skew field}. An \textit{integral domain} is a commutative ring with no nontrivial zero divisors. We will call a domain  \textit{nontrivial} if it is nonzero and not a field. Such domains have zero socles. 

Given a ring $S$ and an $S\text-S$-bimodule $B$, the \textit{trivial extension $S\ltimes B$} is defined as a ring of matrices  \[\{\begin{pmatrix}s&b\\0&s\end{pmatrix}\mid s\in S, b \in B\}.\] 

We begin with the lemma, which provides the local rings used in the constructions
\begin{Lemma}\label{LemmaZ}
    Let $S$ be a nontrivial local integral domain, and let $\m:=S/J(S)$ be its residue skew field. Suppose there exists a ring embedding $\sigma: S\to \m$. On the right $S$-module $\m$, consider the left multiplication $S\times \m\to \m$ by elements $s\in S$ given by $s\cdot m:=\sigma(s)\cdot m$. 
    
    Then $\m$ acquires the structure of an $(S\text-S)$-bimodule, denoted $K:= {_\sigma}\m_S$. Dually,  we consider $(S\text-S)$-bimodule $L:={_S}\m_\sigma$. Furthermore, the following holds:

    (1) As a right module, $K$ is simple, and as a left module, it has a zero socle. 

    (i') As a left module, $L$ is simple, and as a right module, it has a zero socle. 

(2) $N:=K\oplus L$ is an $(S-S)$-bimodule and $soc(_SN)=L$ and $soc(N_S)=K$.

(3) The trivial extension $Z:=S\ltimes N$ is a local ring such that 
$soc({_ZZ})$ and $soc(Z_Z)$ are simple and $soc({}_ZZ)\cap soc(Z_Z)=0$.
\end{Lemma}
\begin{Proof}
    (1) $K_S=\m_S$ is simple because $\m$ is the residue field of $S$. 

   We show that $_SK$ has zero socle by showing that all cyclic submodules have zero socle. Take a nonzero $k\in K$, and consider a nonzero epimorphism $S\twoheadrightarrow Sk$ given by $s\mapsto \sigma(s)k$. Its kernel is the left annihilator of $k$. But for any nonzero $s\in S$, its image $\sigma(s)$  is nonzero in the field $\m$, hence a unit.  Thus, the kernel is trivial, giving $S\cong Sk$. Since $S$ is a nontrivial domain, it has a zero socle, and so does $Sk$.

    (2) It follows from Lemma \ref{LemmaUrelated}.

    (3) A trivial extension of a commutative local ring is local [AW 09, Thm. 3.2(1)]. In a nontrivial local ring, all simple ideals are nilpotent. Because $S$ is a domain, $\left(\begin{smallmatrix}s&n\\0&s\end{smallmatrix}\right)$ is nilpoten  if and only if $s=0$; thus, simple right (left) ideals of $Z$ coincide with simple submodules of a right (left) $S$-module $N$.
\end{Proof}

The following example demonstrates the existence of a nontrivial local integral domain that embeds into its residue field.
\begin{Ex}  Let $\mathbb{X}$ be a countable set of variables, and let $x$ be a variable not in $\mathbb X$. Futher, let $T=\mathbb{Q}(\mathbb X)$ be the fraction field of $\mathbb{Q}[\mathbb X]$. Futhermore, let $S=T[x]_{(x)}$ be a localization of the polynomial ring $T[x]$ in the prime ideal $(x)$. Then $S/xS\cong T$ and  the fraction field $T(x)$ is isomorphic to the field $\mathbb{Q}(\mathbb X\cup \{x\})\cong \mathbb{Q}(\mathbb X)=T$, hence there exist a ring embedding $\sigma:S\to S/xS$.
\end{Ex}

\begin{Ex}\label{ExSocsDifer1}
    Ring $Z$ from Lemma \ref{LemmaZ} has a Nakayama permutation, but $S_r\cap S_l=0$.
\end{Ex}
\begin{Proof}
    For a local ring, having a Nakayama permutation amounts to both $S_r$ and $S_l$ being simple. The rest follows from Lemma \ref{LemmaZ}.
\end{Proof}

We have seen that for rings with coinciding socles, having a Nakayama permutation is equivalent to being a Kasch QF-2 ring. This trivially holds for local rings, but it is not true otherwise. We present a series of examples where we calculate socles of rings based on Lemma \ref{LemmaSimplemodules}.

Let $B, C$ be $(S\text-S)$-bimodules. Then $C$ can be viewed as an $(S\ltimes B\text-S\ltimes B)$-bimodule by viewing it as $0\ltimes C$. Then $(0\ltimes B)C=0=C(0\ltimes B)$ and hence the lattice of right (left) $S\ltimes B$-submodules of $C$ coincides with the lattice of right (left) $S$-submodules of $C$.
\begin{Ex} \label{ExQFKasch}
Let $S, K, L$ as in Lemma \ref{LemmaZ}. Then the the trivial formal matrix ring \[R:=\begin{pmatrix}
    S\ltimes K & L\\
    L & S\ltimes K
\end{pmatrix}~with~
        S_r=\begin{pmatrix}
        0\ltimes soc(K_S) & 0\\
        0 & 0\ltimes soc(K_S)
    \end{pmatrix},\qquad S_l=\begin{pmatrix}
        0&soc({_SL})\\
        soc({_SL}) & 0
    \end{pmatrix}
    \] is a Kasch QF-2 ring, but it does not have a Nakayama permutation.
\end{Ex}

\begin{Ex}\label{ExQFKasch2}
    Let $S, K, L$ as in Lemma \ref{LemmaZ}. Then the trivial matrix ring \[R:=\begin{pmatrix}
        S & K & L\\
        L& S& K\\
        K & L &S
    \end{pmatrix}~with~
        S_r=\begin{pmatrix}
         0 & K & 0\\
        0& 0& K\\
        K & 0 &0
    \end{pmatrix}, \qquad S_l=\begin{pmatrix}
        0 & 0 & L\\
        L& 0& 0\\
        0 & L &0
    \end{pmatrix}\]
  is a QF-2 Kasch ring but it does not have a Nakayama permutation.
\end{Ex}The construction for rings of order $n>3$ is analogous.

As discussed in [Kra. 25, Sec. 4], formal matrix rings with Nakayama permutation and an essential right socle have a well-described structure, generalising the structure of PF rings (Sec \ref{SecMueller}). The following example illustrates the need for assuming the socle's essentiality.
\begin{Ex}\label{ExnotEss}
    Let $S$ be a local commutative ring that is not perfect and $E$ be the injective envelope of its simple module and $M$ be an $S$-module with zero socle.     Then the trivial formal matrix ring \[\begin{pmatrix}
        S&E&M\\
        0&S&E\\
        E\oplus M&0&S
    \end{pmatrix} \]   is a ring with a Nakayama permutation $\pi=(1~2~3)$ and its right and left socles coincide.
\end{Ex}

\section{Double annihilator property}\label{SecDouble}

In this section, we discuss the relationship between various versions of the Double annihilator property and the Nakayama permutation. C. R. Hajarnavis and  N. C. Norton proved that a ring satisfying the Double annihilator property on all one-sided ideals is necessarily a semiperfect ring with a Nakayama permutation and essential socles (Thm. \ref{ThmD}). By contrast, if the ring is semiperfect, it suffices to assume this property only for specific classes of ideals [Nor. 75, Thm. 2.2.9]. The main theorem of this section refines Norton's result. 

We show that a ring with coinciding socles possesses a Nakayama permutation if and only if the Double annihilator property holds for one-sided semisimple ideals or one-sided ideals containing the Jacobson radical. This, in turn, is equivalent to the ring being a Kasch minannihilator ring (Thm. \ref{ThmAnn1}).  We also show that the existence of \textit{any} such lattice anti-isomorphism implies that it is given by annihilators (Thm. \ref{ThmAnti-isom}).

\smallskip 

We first establish properties of ideals with the Double annihilator condition (Sec. \ref{SecDIntro}) and discuss classical results on artinian rings (Sec. \ref{SecDart}). Next, we provide several characterisations of Kasch rings in terms of the Double annihilator property (Sec. \ref{SecDKasch}) and review the literature on minannihilator rings (Sec. \ref{SecDMin}).  Kasch minannihilator rings are then characterised in terms of the Nakayama permutation (Sec. \ref{SecDNak}). The final section  (Sec. \ref{SecDIsom}) then generalises these results for arbitrary lattice anti-isomorphisms.

\subsection{Background and general properties}\label{SecDIntro}

This section introduces annihilator maps $l$ and $r$ and establishes their fundamental property as order-reversing maps between $\mathcal{L}_R(R)$ and $_R\mathcal{L}(R)$. Note that before the notion of Morita duality (Sec. \ref{SecDual}), \textit{the duality of rings} referred to two rings such that their respective lattices of right ideals are anti-isomorphic [Baer 43].

\begin{Def}[D-ideals, D-rings]\label{DefD}
Let $R$ be a ring and $X\subseteq R$ a set. 

We define the \emph{left and right annihilator of $X$}  as \[l(X):=\{r\in R\mid rX=0 \}\quad r(X):=\{r\in R\mid Xr=0\}.\] 

A right (left) ideal $I$ ($I'$) is called a \emph{right (left) D-ideal} if $rl(I)=I$ ($lr(I')=I'$).


A ring is called a \emph{right (left) D-ring} if all right (left) ideals are right (left) D-ideals.  
\end{Def}
We use the name \textit{D-ideals} (shorthand for \textit{dual ideals}), following the convention of [HN 85]. Although the term \textit{annihilator ideal} is quite common, we prefer D-ideal to avoid confusion with the annihilator of a module.
\begin{remark}
    Note the distinction between the annihilator of a set, as defined above, and the concept of the \textit{annihilator of a module}. When the right ideal $I$ is viewed as a right $R$-module $I_R$, its annihilator is $r(I)$, not $l(I)$. In particular, two isomorphic right modules have equal annihilators, whereas the left annihilators of two isomorphic right ideals may differ. 
\end{remark}
\begin{Ex}
    The left annihilator $l(x)$ of an element $x\in R$ coincides with the left annihilator of the right principal ideal $xR$. It is precisely the kernel of the right multiplication homomorphism $\cdot  x\colon R\to R$. In particular, for an idempotent $f$ we have $l(fR)=R(1-f)$.
\end{Ex}
In general, for any subset $X\subseteq R$, the annihilator $l(X)$ is a left ideal and coincides with the left annihilator of the right ideal generated by $X$. 
\begin{Ex}
 The zero ideal and the ring itself are D-ideals.
\end{Ex}
 More examples can be found in Sec. \ref{SecDKasch}, where we interpret the Double annihilator property for various ideals.

The following lemma gathers general properties of maps $l$ and $r$ used implicitly throughout this section.
\begin{Lemma}\label{Lemmalr}
    Let $R$ be a ring. 

 (1)   Then, we have order-reversing maps
\[
\begin{array}{rcl}
_R\mathcal{L}(R) & \overset{r}{\underset{l}{\rightleftarrows}} & \mathcal{L}_R(R) \\
I' & \mapsto & r(I') \\
l(I) & \mapsfrom & I
\end{array}
\]

(2) These maps take finite joins (sums) into meets (intersections): $l(I_1+I_2)=l(I_1)\cap l(I_2)$.

(3) Their compositions $lr$ and $rl$ are closure operators. In particular, for any right ideal $I\subseteq R$, we have $I\subseteq rl(I)$.

(4) A right (left) ideal is a right (left) D-ideal if and only if it is a right (left) annihilator of a left (right) ideal.

(5) A ring $R$ is a D-ring if and only if $l$ and $r$ are both lattice anti-isomorphisms. 
\end{Lemma}
Note that $l$ maps the interval $[J(R), R]_{\mathcal{L}_R(R)}$ to $\mathcal{L}_R(S_r)$ (Lemma \ref{Lemmal(J(R))}) and $r$ maps $\mathcal{L}_R(S_r)$ to $[J(R), R]_{\mathcal{L}_R(R)}$ if the ring is right Kasch rings (Lemma \ref{LemmaAnnSoc}). We explore this connection further in Secs. \ref{SecDNak}, \ref{SecDIsom}.

In contrast to (2), an annihilator of an intersection need not be the sum of the annihilators. However, it is true for right D-ideals such that the sum of their left annihilators is a left D-ideal [HN 85, Lemma 3.1]. Hence, if $l$ and $r$ are surjective, they are anti-homomorphisms. Also, in this case, all ideals are annihilators, so by (4), these maps are injective, which proves (5). 
\begin{remark}
    Rings where maps $l$ and $r$ take finite meets into joins are often called \textit{Ikeda-Nakayama rings} due to [IN 54, Condition (b)]. D-rings are the class of Ikeda-Nakayama rings in which for any simple right $R$-module $T$, the left module $Hom_R(T, R)$ is simple [CNY 00, Thm. 10].  
\end{remark}
\begin{Thm}[HN 85, Thms. 3.5, 3.9 \& 4.2]\label{ThmD}
Let $R$ be any ring.

If $R$ satisfies the Double annihilator property, then it is a semiperfect ring with a Nakayama permutation and essential socles.   
\end{Thm}
Furthermore, all f.g. modules over a D-ring have a finite uniform dimension [HN 85, Thm. 5.3], and a D-ring is perfect if and only if it is a QF rings [HN 85, Cor. 5.5].

\subsection{Artinian case}\label{SecDart}

Nakayama showed that a fin. dim. algebra is  QF if and only if it is a D-ring [Nak. 39, Thms. 1, 2]. Furthermore, it suffices to assume the Double annihilator property only for nilpotent simple ideals and the Jacobson radical [Nak. 39, Thm. 3]. In [Nak. 41a], he observed that this characterisation extends to any artinian ring. While he used the Nakayama permutation to define QF rings, he noted that \textit{\say{this definition does not have much significance for
 itself [in the artinian case], contrary to the case of algebras}} and that this work \textit{\say{may be considered as a study
 of the structure of [artinian D-rings]}.}

 \begin{Thm}\label{ThmDArtin}
Let $R$ be an artinian ring. TFAE

(1) It is a QF ring.

(2) All nilpotent simple right (and left) ideals are right (and left) D-ideals, and $J(R)$ is a D-ideal.

(3) All right (and left) ideals are right (and left) D-ideals.

(4) There are anti-isomorphisms \[ {}_R\mathcal{L}(top(R)) \cong_{op} \mathcal{L}_R(S_r)\quad  {}_R\mathcal{L}(S_l) \cong_{op} \mathcal{L}_R(top(R)).\]
\end{Thm}
\begin{Proof}
    The equivalence (1)$\Leftrightarrow$(2)$\Leftrightarrow(3)$ is [Nak. 41a, Thm. 6]. The equivalence (1)$\Leftrightarrow$(4) is [Nak. 41b, Thm. 1].
\end{Proof}

By Lemma \ref{Lemmalr}, condition (3) induces an anti-isomorphism \[_R\mathcal{L}(R) \cong_{op} \mathcal{L}_R(R),\]  so it suffices to assume that $R$ is right artinian. Furthermore, the dualities in (4) are thus given by the annihilator maps. By (1), it follows that $S_l=S_r$ (Prop. \ref{PropQFKaschEss})

Note that a simple right ideal is either nilpotent or generated by an idempotent. Hence, in view of Thm. \ref{ThmKasch}(4),  condition (2) states that the ring is a Kasch minannihilator ring.
\begin{Def} Let $R$ be a ring.

Then $R$ is a \emph{right (left) minannihilator ring} if all simple right (left) ideals are right (left) D-ideals. 
\end{Def}
Our main result (Thm. \ref{ThmAnn1}) is a generalisation of the equivalence  (1)$\Leftrightarrow$(2). Although Nakayama's proof [Nak. 39, Thm. 3] is also valid for semiprimary rings, further generalisation requires different methods, as his approach relies on the nilpotency of the Jacobson radical. In the artinian case, it suffices to assume that the ring is a minannihilator ring (Prop. \ref{PropStorrer}). We generalise the equivalence (1)$\Leftrightarrow$(4) in Thm. \ref{ThmAnti-isom}.

\subsection{Kasch rings}\label{SecDKasch}

Recall that a ring $R$ is \textit{right Kasch} if and only if any simple right module $T$ has an isomorphic copy in the right socle of $R$.  This section discusses various known characterisations of semiperfect Kasch rings in terms of annihilators, culminating in the following theorem:

\begin{Thm}\label{ThmKasch}Let $R$ be a semiperfect ring. TFAE

(1) $R$ is a right Kasch ring. 

(2) All maximal right ideals are right D-ideals.

(3) All right ideals containing $J(R)$ are right D-ideals.

(4) $J(R)$ is a right D-ideal.     
\end{Thm}
Artinian Kasch rings were introduced under the name \textit{S-rings} by Kasch in [Kas. 54] while establishing the basics of the
theory of Frobenius extensions. 

We begin by citing the following basic properties gathered by Morita, noting that some of these were already known.
\begin{Prop}[Mor. 66, Thm. 1]\label{PropKasch}
    Let $R$ be a ring and $M\subseteq R$ a maximal right ideal. TFAE

    (1) $M$ is a right D-ideal.

    (2) $l(M)\neq 0$.

    (3) $R/M$ has an isomorphic copy inside $soc(R_R)$.

    (4) $Hom_R(R/M,  R)$ is nonzero.
\end{Prop}
Condition (2) explains why \textit{right} Kasch rings used to be called \textit{left} S-rings.

 For semilocal rings, condition (1) can be further strengthened. Norton notes that the following is well-known
\begin{Lemma}[Nor. 75, Lemma 2.2.6]\label{LemmaAboveJ}
    Let $R$ be a semilocal right Kasch ring.

 Then, for any right ideal $I$, we have $rl(I)+J(R)=I+J(R)$. In particular, right ideals containing $J(R)$ are right D-ideals.
    \end{Lemma}
  Thus, the Jacobson radical of a right Kasch ring is a right D-ideal. The converse is shown in Lemma \ref{LemmaAnnSoc}.
    
    In general, the annihilator of a maximal right ideal does not need to be simple (cf. Lemma \ref{LemmaBijectionSimple}), as discussed by the following lemma. 
    
    Recall that in a semiperfect ring, all right ideals containing the radical can be written as $(1-f)R+J(R)$ for some idempotent $f\in R$ (Lemma \ref{LemmaIdealsAboveJ}).

     \begin{Lemma}\label{LemmaAnnMax}
    Let $R$ be a semiperfect ring and $I=(1-f)R+J(R)$ a right ideal.

    If $S_l=S_r$, then $l(I)=soc(Rf)$.
\end{Lemma}
\begin{Proof}
We have \[l(I)=l((1-f)R+J(R))=Rf\cap S_r=Rf\cap S_l=soc(Rf),\]
using Lemma \ref{LemmaSoceR}.
\end{Proof}

Before proceeding further, we recall the following well-known property of the Jacobson radical.
\begin{Lemma}[Lam 91, Cor. 4.2]\label{LemmaJdef}
    Let $R$ be a ring.

    Then the Jacobson radical is the set of elements $r$ such that $Mr=0$ for all simple right $R$-modules $M$.
    
        Furthermore, it is enough to test this property only on a set of representatives of the isomorphism classes of simple right modules.
\end{Lemma}

This allows us to describe the annihilator of the Jacobson radical.
\begin{Lemma}\label{Lemmal(J(R))}
    If $R$ is semilocal, then  $soc(R_R)=l(J(R))$.
\end{Lemma}
\begin{Proof}
Because $J(R)$ is an ideal, $l(J(R))$ is a right ideal. It is annihilated by $J(R)$, so it has the structure of a right $R/J(R)$-module, where $R/J(R)=top(R)$. Thus $l(J(R))$ is semisimple as a right $R$-module since $R$ is semilocal.

On the other hand,  $TJ(R)=0$ for any simple right module $T$, so the right socle is contained in $l(J(R))$.   
\end{Proof}
   
We finish by describing the right annihilator of the right socle.
\begin{Lemma}\label{LemmaAnnSoc}
If $R$ is right Kasch, then $r(soc(R_R))=J(R)$.

If $R$ is semiperfect, then it is right Kasch if and only if $r(S_r)=J(R)$.
\end{Lemma}
The second part of the lemma is adapted from [Die. 58, (4.2)].
\begin{Proof}
 An element $a\in r(soc(R_R))$ annihilates all simple right modules that inject into the right socle. Thus, if $R$ is right Kasch, then $a$ annihilates all right modules, proving $a\in J(R)$ by Lemma \ref{LemmaJdef}. The same lemma then also shows the opposite inclusion. 

Now assume that the ring $R$ is semiperfect and $r(S_r)=J(R)$. Consider the Wedderburn-Artin decomposition of $top(R)$ and let $\bar{e_k}$ be the unities corresponding to each ring $M_{\mu_k}(\m_k)$ of this decomposition. 

Each $\bar{e_k}$ is an idempotent, so they lift modulo the radical in $R$, yielding a decomposition $1_R=\sum e_k$, a sum of pairwise orthogonal idempotents $e_k$, which are the preimages of $\bar{e_k}$.

Nonzero idempotents are not contained in $J(R)=r(S_r)$, so each $S_re_k$ is nonzero. Now $S_re_k$ is the homogeneous component $S_k$ of $S_r$. Indeed, any simple right ideal $T$ of $R$ is also a right $top(R)$-module. So $e_k$ annihilates all simple right ideals which are not isomorphic to $V_k$. 

In total, each homogeneous component $S_k$ of $S_r$ is nonzero, so the ring is right Kasch. 
\end{Proof}

\subsection{Minannihilator rings}\label{SecDMin}

There is a rich literature on minannihilator rings. We gather results concerning the Nakayama permutation. We also introduce the concept of \textit{mininjective rings}, which generalises several classical results.  Results of this section will not be used in the proof of Theorem \ref{ThmAnn1}, but it is worth noting that a crucial step in the proof of Theorem \ref{ThmNY} by W.  Nicholson and M. Yousif is showing that such rings are Kasch.

\smallskip  

Minannihilator rings were first studied in the artinian context. M. Ikeda proved that they are QF rings [Ike. 51, Cor. 3]. Storrer further refined this result:
\begin{Prop}[Sto. 69, Prop. 1]\label{PropStorrer} 
Let $R$ be a right artinian ring.

Then $R$ is QF if and only if it is a minannihilator ring.     
\end{Prop}
Storrer's proof relies not on the minannihilator property, but on the following consequence (cf. Prop. \ref{PropKasch}(4)): 
\begin{Lemma}
    Let $R$ be a right minannihilator ring with an essential left socle and let $T$ be a simple left $R$-module.

    Then the right $R$-module $Hom_R(T, R)$ is either simple or zero.
\end{Lemma}
This lemma was originally formulated by Dieudonné [Die. 58, 3.2] for artinian rings, but the proof works in this more general setting. 

\smallskip 

minannihilator rings are closely related to \textit{mininjective} rings. Note that a left minannihilator ring with an essential left socle is right mininjective [NY 97, Cor. 2.6].
\begin{Def} Let $R$ be a ring.

Then $R$ is \emph{right mininjective} if for any simple right ideal $T$ in $R$, any $R$-homomorphism $T\to R$ extends to a homomorphism $R\to R$.
\end{Def}
Historically, the concept of a mininjective module was introduced by Harada in [Har. 82]where a right mininjective ring is referred to as a \textit{right self mini-injective ring}. For rings, this condition was introduced by Ikeda in [Ike. 52, Lemma 1], as a special case of the \textit{Shoda condition}, and he showed that a right mininjective fin. dim. algebra over a field is QF [Ike. 52, Prop. 1]. The following theorem, due to Nicholson and Yousif, strengthens Ikeda's observation:
\begin{Thm}[NY 03, Thm. 3.31]
    Let $R$ be a semilocal mininjective ring
with the ascending chain condition on right annihilators and an essential right socle. 

Then $R$ is QF.
\end{Thm}
For a discussion on newer improvements to this theorem, see [NKC 23].

The following proposition explains the connection between mininjective rings and minannihilator rings:
\begin{Prop}[NY 97, Prop. 2.4]\label{PropNY}  Let $R$ be a right minannihilator ring.

Then $R$ is left mininjective if and only if the left socle is contained in the right socle.    
\end{Prop}

This leads to the central result of [NY 97], which describes semiperfect mininjective rings.
\begin{Thm}[NY 97, Cor. 3.10]\label{ThmNY}
Let $R$ be a semiperfect mininjective ring. Further, suppose that all primitive right and all primitive left ideals have a nonzero socle.

Then $R$ is a ring with a Nakayama permutation and coinciding socles.    
\end{Thm}

In particular, the assumption on primitive ideals in the previous theorem is satisfied if both socles are essential. Together with Prop. \ref{PropNY} this yields:
\begin{Cor}
    Let $R$ be a semiperfect minannihilator ring with coinciding essential socles. 

    Then $R$ has a Nakayama permutation.
\end{Cor}

\subsection{Nakayama permutation}\label{SecDNak}

A D-ring is a semiperfect ring with a Nakayama permutation and essential socles. For a semiperfect ring, it is enough to assume that only \textit{some} one-sided ideals are annihilators. This is implicit in [HN 85, Secs. 3, 4]. The converse is not true; Norton provides an example of a commutative local ring with essential socles and a Nakayama permutation, which is not a D-ring [Nor. 75, Ex. 3.2.3].

The main theorem of this section (Thm. \ref{ThmAnn1}) characterises Kasch minannihilator rings as rings with a Nakayama permutation and coinciding socles. As a corollary, we show that in such rings, all semisimple right ideals are of the form $fsoc(R)$ for some idempotent $f$ (Cor. \ref{CorDmain}).

We use the methods of [Nor. 75] to prove the implication (2)$\implies$(1) of Thm. \ref{ThmAnn1}.  As far as we are aware, implication (1)$\implies$(3) is new.

\smallskip 

We begin by recalling the following Lemma:
\begin{Lemma}\label{LemmaBijectionSimple}
    Let $R$ be a ring.
    
      Then $R$ is a left Kasch right minannihilator ring if and only if $l$ and $r$ are mutually inverse bijections between the sets of simple right and maximal left ideals.
\end{Lemma}
\begin{Proof}
If such bijections exist, all left maximal (right simple) ideals are left (right) D-ideals, and hence $R$ is a left Kasch (right minannihilator) ring.

For the converse,  let $T$ be a right simple ideal and $M$ a maximal left ideal such that $l(T)\subseteq M$. Applying $r$, we get $T=rl(T)\supseteq r(M)$. Thus $r(M)$ is contained in a simple ideal. Since it is nonzero by the left Kasch property (Prop. \ref{PropKasch}(3)), we conclude $r(M)=T$. Applying $l$ again, we get that $M=lr(M)=l(T)$.
\end{Proof}

Due to symmetry, we obtain a dual result by switching the domains of the maps $l$ and $r$. It turns out that assuming both of these versions leads to a stronger result. 
\begin{Thm}\label{ThmAnn1} Let $R$ be a semiperfect ring. TFAE

(1)   $R$ possesses a Nakayama permutation and $S_l=S_r$.

(2) $R$ is a Kasch minannihilator ring.

(3) All right (left) ideals contained in the right (left) socle and all right (left) ideals containing the radical are right (left) D-ideals.    

(4) The maps $l$ and $r$ induce mutually inverse lattice anti-isomorphisms:
\[
    {}_R\mathcal{L}(top(R)) \cong_{op} \mathcal{L}_R(S_r)
    \quad \text{and} \quad
  {}_R\mathcal{L}(S_l)  \cong_{op} \mathcal{L}_R(top(R)).\]
\end{Thm}
Applying Lemma \ref{LemmaAnnMax}(4) gives the following corollary.
\begin{Cor} \label{CorDmain} Let $R$ be a semiperfect ring with a Nakayama permutation and $S_r=S_l$, and let $I$ be a semisimple right ideal.

Then there exists an idempotent $f$ such that $I=soc(fR)=fsoc(R)$.    
\end{Cor}
We now extend Lemma \ref{LemmaBijectionSimple}.
\begin{Lemma}\label{LemmaBijectionSemiSimple}
    Let $R$ be a semiperfect left Kasch ring such that all semisimple right ideals are right D-ideals. 

    If $S_r=S_l$, then $l$ and $r$ induce mutually inverse lattice anti-isomorphisms:
\[{}_R\mathcal{L}(top(R) \cong_{op} \mathcal{L}_R(soc(R))\]
\end{Lemma}
Recall that the all left ideals containing $J(R)$ are left D-ideals if and only if the ring is left Kasch (Thm. \ref{ThmKasch}).
\begin{Proof}  If $I'\supseteq J(R)$ is a left ideal, then $r(I')\subseteq r(J(R))=soc(R)$.   Let $I\subseteq soc(R)$ be a right ideal. Then $l(I)\supseteq l(soc(R))$. Because $R$ is left Kasch, we have $l(soc(R))=J(R)$ (Lemma \ref{LemmaAnnSoc}).   So $l$ and $r$ are mutually inverse bijections 
\[{}_R\mathcal{L}(top(R)) \underset{l}{\overset{r}{\rightleftarrows}} \mathcal{L}_R(soc(R)).\]
Because they are order-reversing maps, they are lattice anti-isomorphisms. The conclusion then follows from the identification 
\[ [J(R), R]_{_R\mathcal{L}(R)} = {}_R\mathcal{L}(top(R)).\]
\end{Proof}

\smallskip

The structure of the proof of Thm. \ref{ThmAnn1} is as follows: Implication (4)$\implies$(3) is trivial. Once we prove that (3)$\implies$(2)$\implies$(1)$\implies$(3), the implication (3)$\implies$(4) follows directly from the previous lemma. Implication (3)$\implies$(2)  follows from the definitions. 

So it remains to prove (2)$\implies$(1)$\implies$(3).  We first prove (2)$\implies$(1). First, we modify  [HN 85, Lemma 3.1,  3.2].
\begin{Lemma}
Let $R$ be a right Kasch ring. 

If $J$ is a superfluous right D-ideal, then $l(J)$ has a nonzero intersection with all nonzero left D-ideals. 
\end{Lemma}
\begin{Proof}
    Let $X$ be a left D-ideal and consider $l(J)\cap X$. Applying $r$ we get
    \[r(l(J)\cap X)=r(lrl(J)\cap lr(X))=rl(rl(J)+r(X))=rl(J+r(X)),\]
where the first equality follows because $X$ and $l(J)$ are left D-ideals. 

If $l(J)\cap X=0$, we have $R=r(l(J+r(X)))$ and hence $0=l(J+r(X))$. Since $R$ is right Kasch, the only ideal with a trivial left annihilator is $R$ itself, i.e., $R=J+r(X)$. But $J$ is assumed to be superfluous, hence $r(X)=R$, proving that $X=0$.
\end{Proof}
\begin{remark}
    In particular, the Jacobson radical is a superfluous ideal. Hence, by Lemma \ref{Lemmal(J(R))}, if $R$ is semilocal, we know that the right socle has a nontrivial intersection with left D-ideals. This shows that the conditions in Theorem \ref{ThmAnn1} are not equivalent to $R$ being a D-ring, because such rings necessarily have essential socles; thus, Example \ref{ExnotEss} serves as a counterexample.
\end{remark}

We can now continue the proof of the main theorem. 
\begin{Proof}[(2)$\implies$(1)]
Because all simple left ideals are assumed to be left D-ideals, and the ring is right Kasch, by the previous lemma, $S_l\subseteq l(J)=S_r$. Similarly, $S_r\subseteq r(J)=S_l$.

Thus $R$ has coinciding socles. Now for an idempotent $e\in R$, we have $l((1-e)R+J(R))=soc(eR)$ by Lemma \ref{LemmaAnnMax}. If $e$ is primitive, then $soc(eR)$ is the annihilator of a maximal right ideal. By Lemma \ref{LemmaBijectionSimple}, this means that it is a simple right ideal, proving that $R$ is a right QF-2 ring. 

    In an analogous manner, we prove that it is left QF-2. Hence, we have a QF-2 ring with coinciding socles, proving that it has a Nakayama permutation by Proposition \ref{PropSocsCoincide}.
\end{Proof}

We now aim to prove that (1)$\implies$(3). We will repeatedly use the fact, following from the definition of a ring with a Nakayama permutation, that such rings are Kasch. By Lemma \ref{LemmaAboveJ}, right (left) ideals containing $J(R)$ are right (left) D-ideals; therefore, it is sufficient to prove that right (left) semisimple ideals are right (left) D-ideals. 

\begin{Lemma}\label{LemmaAnnsurj}
    Let $R$ be a semiperfect left Kasch ring such that $S_r=S_l$. 

    Then $l$ is a surjection 
    \[\mathcal{L}_R(soc(R)) \overset{l}{\twoheadrightarrow} [J(R), R]_{_R\mathcal{L}}.\]
\end{Lemma}
\begin{Proof}
Because $R$ is left Kasch, $l(I)$ contains $J(R)$ for any semisimple right ideal $I$ (Lemma \ref{LemmaAnnSoc}). To show that $l$ is surjective, recall that in a left Kasch ring, any left ideal $I'\supseteq J(R)$ can be written as $I'=l(I)$ for some right ideal $I$ in  $R$ (Lemma \ref{LemmaAboveJ}). This $I$ is semisimple as $soc(R)=r(J(R))\supseteq r(I')=rl(I)\supseteq I$. 
\end{Proof}

It remains to prove that $l$ is also an injection. We recall that in our notation, the composition length of the module $top(R)$ equals $m$, and if the ring possesses a Nakayama permutation, then the composition length of the right socle is $m$ too.

\begin{Proof}[(1)$\implies$(3)] We prove any pair of distinct right semisimple ideals has distinct left annihilators. It is sufficient to consider the case where one ideal is contained in the other. Indeed, let $I, J\subseteq soc(R)$ be right ideals such that $l(I)=l(J)$. If neither of them is contained in the other, then $I$ is properly contained in $I+J\subseteq soc(R)$. But $l(I+J)=l(I)\cap l(J)=l(I)$.

Without loss of generality, consider two right ideals $I\subseteq J\subseteq soc(R)$ such that $J/I$ is simple. We can extend $I\subseteq J$ to a composition series
\[0=I_0\subseteq I_1\subseteq \dots \subseteq I_\lambda=I\subseteq J=I_{\lambda+1}\subseteq \dots \subseteq I_{m}=soc(R),\]
    for some $\lambda\subseteq m$.

   Applying $l$ gives a series 
\[R=l(I_0)\supseteq l(I_1)\supseteq \dots \supseteq l(soc(R)),\]
where consecutive factors are either simple or zero. 

Because $R$ has a Nakayama permutation,  $top(R)$  has the same composition length as $soc(R)$, and the ring is Kasch, showing $l(soc(R))=J(R)$ (Lemma \ref{LemmaAnnSoc}). Hence, this series is a composition series, and none of the consecutive factors is zero. In particular, $l(J)$ is properly contained in $l(I)$.
\end{Proof}

\subsection{Anti-isomorphisms of lattices}\label{SecDIsom}

This section shows that a semiperfect ring with essential socles possesses a Nakayama permutation if and only if the lattices of submodules of socle and top are anti-isomorphic (Thm. \ref{ThmAnti-isom}), generalising Nakayama's result (Thm. \ref{ThmDArtin}(4)). Thus, for a semiperfect ring with essential socles, if such an anti-isomorphism exists, it is given by annihilator maps. 

For this section, recall $m=\mu_1+\dots +\mu_n$ denotes the maximal number of pairwise orthogonal primitive idempotents in the decomposition of the unity. Also recall the notation from Section \ref{Secsemptop}. The proof of Thm. \ref{ThmAnti-isom} is based on Lemmas  \ref{LemmaBir}, \ref{LemmaSocleStruct}.

\smallskip

\begin{Thm}\label{ThmAnti-isom}
Let $R$ be a semiperfect ring with essential socles. 

Then $R$ has a Nakayama permutation if and only if there are anti-isomorphisms  \[ {}_R\mathcal{L}(top(R)) \cong_{op} \mathcal{L}_R(S_r)\quad  {}_R\mathcal{L}(S_l) \cong_{op} \mathcal{L}_R(top(R)).\]
\end{Thm}
\begin{Proof}
First, assume $R$ has a Nakayama permutation and essential socles. Then its left and right socle coincide (Thm. \ref{PropQFKaschEss}) and the conclusion follows from Thm. \ref{ThmAnn1}(4).

We now prove the converse. We use ${}_R\mathcal{L}(top(R)) \cong_{op} \mathcal{L}_R(S_r)$ to show that the ring is right QF-2 and right Kasch. The left version is proved analogously using the second anti-isomorphism. 

Recall that $_R\mathcal{L}(top(R))$ is a modular lattice of finite length. By Lemma \ref{LemmaBir}(4) its rank equals its composition length, namely $m$. Hence, by the assumption, the rank of $\mathcal{L}_R(S_r)$ is also $m$. Thus the composition length of $S_r$ is $m$. By Lemma \ref{LemmaSoceR},  the composition length of $S_r$ is the sum of the composition lengths of the socles of right primitive ideals.  Since $R$ has an essential right socle, each $e_iR$ has a nonzero socle, contributing at least $1$ to the total length. Since the sum of these lengths equals $m$, each $soc(e_iR)$ must have length exactly $1$, i.e., it is simple. We have shown that $R$ is right QF-2. 

In view of Thm. \ref{PropQFKaschEss}, it remains to show that the ring is Kasch. By Lemma \ref{LemmaBir}, the lattice $_R\mathcal{L}(top(R))$ decomposes as a product of $n$ indecomposable lattices, and so does $\mathcal{L}_R(S_r)$  by the assumption. Thus, the socle must contain non-zero homogeneous components corresponding to all $n$ isomorphism types (Lemma \ref{LemmaSocleStruct}).
\end{Proof}
\begin{remark}
    The existence of these lattice anti-isomorphisms is structurally consistent with Prop. \ref{PropIsomfields}. While the lattices of second powers of simple modules are isomorphic based solely on the cardinality of their base fields, a stable correspondence that holds for arbitrary multiplicities $\mu_k$ requires the underlying skew fields $\m_k$ and $\m_{\pi(k)}$ to be isomorphic by the Fundamental Theorem of Projective Geometry [Baer 52, Chap. III.I].
\end{remark}

\section{Cardinality conditions}\label{SecCard}

This section investigates two versions of a \textit{cardinality condition} and their connection to Frobenius rings, which we generalise as \textit{semiperfect rings with a Nakayama permutation that respects the multiplicities}. Cardinality conditions imply the Double annihilator property (Lemma \ref{LemmaDual}). We have seen that the Double annihilator property encodes that the socle and the top have anti-isomorphic lattices of submodules (Thm. \ref{ThmAnti-isom}). We now show that cardinality conditions encode that the socle and the top are isomorphic as modules (Cor. \ref{CorCardmain}).

For a finite ring $R$, one of the cardinality conditions amounts to studying right ideals $I$ and left ideals $I'$ such that \[|l(I)||I|=|R|\qquad |I'||r(I')|=|R|.\] Following  [Wood 09, Sec. 12.4], we say that such $I$ ($I'$) satisfies the \textit{Size condition}. The second version, which we introduce in Definition \ref{DefGenDim}, originates from Nakayama's functions $d_r$ and $d_l$ [Nak. 41a, page 8] and will be referred to as the \textit{Generalised dimension condition}. 

There is a significant overlap between these conditions, as discussed in the sequel, and many results have proofs that are nearly identical, though one is multiplicative and the other is additive. Thus, when a result holds for both of them, we use the term \textit{cardinality condition}, noting that the Size condition only applies to finite rings.

The history of the Size condition dates back to the work of E. Lamprecht [Lam. 53]. Honold used Lamprecht's results on character modules in his characterisation of finite Frobenius rings [Hon. 01, Thm. 1]. He points out that \say{\textit{[Lamprecht's] papers are apparently not as well-known
among ring theorists as they deserve to be.}} However, as the concept of a \textit{character of an abelian group} is well-known, character modules have reappeared independently in ring theory [CG 92, Sec. 3] and were subsequently used to characterise finite Frobenius rings [Hir. 97, Thm. 1].

\smallskip

Section \ref{SecCardDef} conceptualises Nakayama's function $d_r$ as a length function on an abelian category of finite-length modules over a semiperfect ring. We then explain that over finite QF rings, Nakayama's cardinality condition coincides with the Size condition. Section \ref{SecCardClass}  gathers classical results.  We then further investigate cardinality conditions for rings that are assumed to have the Double annihilator property (Sec. \ref{SecCardWood}) and general rings (Sec. \ref{SecCardHon}).

\smallskip

Note that a cardinality condition implies the Double annihilator property.  Due to repeated use of this observation, we formulate it as a lemma. 
\begin{Lemma}\label{LemmaDual}
    If a right ideal $I$ and its double dual $rl(I)$ satisfy a cardinality condition, then $I=rl(I)$.
\end{Lemma}
We prove the version for the Size condition as an illustration.
\begin{Proof}
    Let $I$ be a right ideal and $J$ its double dual. We have $|l(I)||I|=|R|=|l(J)||J|=|R|$. Because $I\subseteq J$ and $l(I)=lrl(I)=l(J)$, we conclude that $|I|=|J|$ and hence $I=J$.
\end{Proof}
In particular, if $M$ is a maximal right ideal, then the cardinality condition implies $rl(M)=M$.

\subsection{Definition and basic properties}\label{SecCardDef}

The origin of the dimension condition goes back to [Nak. 39]. The subsequent paper [Nak. 41a] presented, albeit in the artinian case, a possible way to extend the notion of dimension to an arbitrary semilocal ring by a function $d_r$ (also called \textit{colength} in [Die. 58, Sec. 5]). We define this function (which we simply call $d$) as a length function (in the sense of [Krs. 07]) and motivate and define a \textit{Generalised dimension condition} for rings that are not of finite length, focusing rather on their semisimple ideals.

\smallskip 

Suppose that a finite ring $R$ is a fin. dim. algebra over a field $K$ and let $I$ in $R$ be a right ideal. The Size condition amounts to  $|l(I)||I|=|R|$  which holds if and only if the following  \textit{Dimension condition} holds
\[dim_K(l(I))+dim_K(I)=dim_K(R),\]
as shown by observation that $|I|=|K|^{dim_K(I)}$. 

Naturally, one can consider the Dimension condition even when $K$ is not finite. Furthermore, using a generalisation of the Rank-nullity theorem, it can be restated as 
\[dim_K(I)=dim_K(R)-dim_K(l(I))=dim_K(R/l(I)).\]

This allows us to further extend the Dimension condition to infinite-dimensional algebras as long as we only consider right ideals $I$ such that both $dim_K(I)$ and $dim_K(R/l(I))$ are finite. Recall that in a D-ring,  lattices $\mathcal{L}_R(I)$ and $_R\mathcal{L}(R/l(I))$ are anti-isomorphic (Lemma \ref{Lemmalr}). So, whenever $dim_K(I)$ is finite, $dim_K(R/l(I))$ is also finite.  In view of Theorem \ref{ThmAnn1}, for a ring with a Nakayama permutation and coinciding socles, we can consider a cardinality condition on the right semisimple ideals.

Following [Krs. 07, Sec. 3],  a \textit{length function} on the category of finite-length right (left) $R$-modules can be defined by specifying the lengths on (the representative set of isomorphism classes of) simple right (left) modules and extending it additively to an arbitrary module of finite composition length via the Rank–nullity theorem. That is, for a right (left) ideal $I$ of composition length $\lambda$ with a composition series
\[   0=I_0\subseteq I_1\subseteq \dots\subseteq  I_{\lambda-1}\subseteq I_\lambda=I,\]
 \textit{the length of $I$} is defined as the sum of the lengths of the simple modules $I_i/I_{i-1}$ for $1\leq i\leq \lambda$. By the Jordan-Hölder Theorem, this value is independent of the choice of the composition series.

\begin{Ex}
    The \textit{composition length} is the length function induced by assigning $1$ to all simple modules. 
    \end{Ex}
    For a semiperfect ring $R$, recall the representative set of simple modules  $V_1, \dots, V_n$ with multiplicities $\mu_1,\dots, \mu_n$ (Sec. \ref{SecSempNot}).
    \begin{Ex}
        If $R$ is a $K$-algebra, assigning the $K$-dimension to each simple module yields the $K$-dimension as the induced length function. Note that $dim_K(V_k)=\mu_k dim_K(\m_k)$ for each $k\leq n$. If $K$ is algebraically closed, then $\m_k\cong K$.
    \end{Ex}
\begin{remark}
       The \textit{size} of a finite module is not a length function, as $|M|=|K||L|$ for a short exact sequence $0\to K\to M\to L\to 0$.
\end{remark}

\begin{Def}\label{DefGenDim} [Generalised dimension]
    Let $R$ be a semiperfect ring.

    Then the \emph{generalised dimension}, denoted by $d$, is defined as a length function on the category of finite-length right $R$-modules induced by assigning $\mu_k$ to each simple right $R$-module $V_k$.

    We then say that the right ideal $I\subseteq R$  satisfies the \emph{Generalised dimension condition} if $d(I)$ and $d(R/l(I))$ are finite and equal.
    \end{Def}
    The left versions are defined analogously. 
\begin{Ex}\label{ExdimVk}
Recall that, as an $top(R)$-module, the simple right (left) $R$-module $V_k$ ($V'_k$) corresponds to a row (column) ideal in the ring $M_{\mu_k}(\m_k)$. So $\m_k$ is its ring of endomorphisms, and $dim_{\m_k}(V_k)=\mu_k=dim_{\m_k}(V'_k)$. If $\m_k$ is a finite field, then $|V_k|=|\m_k|^{\mu_k}$.
\end{Ex}
Thus, our generalised dimension coincides with Nakayama's functions $d_r$ and $d_l$. If $R$ is a $K$-algebra, it coincides with the $K$-dimension of a module. 

We can define a cardinality condition with respect to any length function, but it seems that, apart from the generalised dimension, there is no known application for it.
\begin{Ex}
     The cardinality condition with respect to composition length is equivalent to the Double annihilator property.
\end{Ex}

\subsection{Classical results} \label{SecCardClass}

For understanding the role of the Generalised dimension condition, the following well-known observation is key.
\begin{Prop}[Lam. 99, Prop. 16.14]\label{PropLamQF}
Let $R$ be an artinian ring. 

Then $R$ is Frobenius if and only if $S_l\cong top({_R}R)$ and $S_r\cong top(R_R)$. Furthermore, if $R$ is QF, it is enough to assume only one of the isomorphisms.    
\end{Prop}
The isomorphisms automatically imply that the ring is Kasch, and by a counting argument (similar to arguments we are about to see), we then find that socles of primitive ideals are simple. The \say{Futhermore} part simply restates the observation that Frobenius rings are QF rings whose Nakayama permutation respects multiplicities.

This proposition can be sharpened in the case of fin. dim. algebras. 
\begin{Prop}[Nak. 49, Cor. 2] \label{PropNakPrinc}
Let $R$ be a finite-dimensional $K$-algebra.

Then $R$ is a Frobenius ring if and only if $S_r$ is a principal right ideal.     
\end{Prop}
\begin{remark}
    If $R$ is a semiperfect ring whose all primitive right ideals have a nonzero socle, then the condition $S_r\cong top(R_R)$ is equivalent to $S_r$ being principal. Indeed, a (semisimple) right ideal is a cyclic module and thus an epimorphic image of (the top of) $R_R$. Hence, $S_r$ being principal implies that its composition length is at most $m$, but because all primitive right ideals have nonzero socles, it is also at least $m$. Thus, $S_r$ is an epimorphic image of $top(R_R)$, but it has the same finite composition length, and thus they are isomorphic.
\end{remark}

Nakayama follows [Nak. 49, Cor. 2] by providing an example that shows it is not true for a general artinian ring. In the late 1990s, Xue proved the version of the above proposition for finite local rings [Xue 98]. A few years later, Honold generalised it to any finite ring [Hon. 01,  Thm. 1]. Finally, in the last decade, Iovanov gave a common generalisation of these propositions.
\begin{Thm}[Iov. 16, Thm. 1.7]\label{ThmIovCyclic}
Let $R$ be an artin algebra.

Then $R$ is Frobenius if and only if $S_r\cong top(R_R)$ as right $R$-modules.    
\end{Thm}
Iovanov's proof uses the existence of a Morita self-duality for Artin algebras, pointing out why this fails for a general artinian ring.

We now present the initial result by Nakayama:
\begin{Thm}[Nak. 41a, Thm. 7] Let $R$ be an artinian ring. TFAE\label{ThmNakCard}

(1) $R$ is Frobenius. 

(2) $R$ is QF and the Generalised dimension condition holds for all nilpotent simple right and nilpotent simple left ideals.

(3) The Generalised dimension condition holds for all right and all left ideals.    
\end{Thm}
 Honold studied the Size condition without assuming the Double annihilator property:
\begin{Thm}[Hon. 01, Thms. 1, 2]\label{ThmHonCard}
Let $R$ be a finite ring. TFAE

(1) $R$ is Frobenius. 

(2) The Size condition holds for all right and all left ideals.        

(3) The Size condition holds for all right ideals.  

(4) The Size condition holds for the Jacobson radical and homogeneous components of the right socle viewed as right ideals, and for maximal right ideals.
\end{Thm}
 Nakayama's theorem does not directly imply Honold's, as the former is dimension-based while the latter is cardinality-based. So one needs to have control over the sizes of skew fields $\m_k$ by means of Prop. \ref{PropIsomfields}.

\subsection{Woods' criterion, its interpretation and generalisation}\label{SecCardWood}

This section examines cardinality conditions on rings, assuming the Double annihilator property for semisimple ideals. In view of Thm. \ref{ThmAnn1}, these are precisely rings with a Nakayama permutation and coinciding socles.  While the Double annihilator property is tied to homological properties of the ring, the cardinality condition is tied to its combinatorics, which can differ even among Morita equivalent rings, as it takes into account changes in multiplicities. The cardinality condition then strengthens the existence of a Nakayama permutation by imposing that the multiplicities in a cycle of the Nakayama permutation coincide.  

\smallskip

First, we generalise the concept of a Frobenius ring. 

\begin{Def}\label{DefMult}
    Let $R$ be a semiperfect ring with a Nakayama permutation $\pi$.
    
    We say that $\pi$ \emph{respects the multiplicities of $R$} if for any $k\leq n$ the equality $\mu_k=\mu_{\pi(k)}$ holds.
\end{Def}
 A corner $eRe$ of $R$ induced by a cycle in $\pi$ has a cyclic Nakayama permutation (Thm. \ref{ThmCycle}). Respecting multiplicities over this cycle means that $eRe$ is a (classical) matrix ring over its basic ring.

\smallskip 

 Wood proved the following proposition, providing an alternative proof that finite rings satisfying the Size condition are Frobenius.

\begin{Thm}[Wood 09, Thm. 12.1]\label{ThmWoodSimple}
Let $R$ be a finite QF ring that is not Frobenius. 

Then there exists a simple right ideal $T$ and a simple left ideal $T'$ such that $|l(T)||T|>|R|$ and $|T'||r(T')|>|R|$.   
\end{Thm}

Thus, if we know that a ring is QF, it is enough to check the Size condition on the simple right ideals (cf. Prop. \ref{ThmNakCard}(2)). Recalling Lemma \ref{LemmaDual}, it follows that rings with the Size condition are Frobenius (cf. Prop. \ref{ThmNakCard}(3)). 

However, if we do not assume that the ring is a QF ring, then the Size condition on simple one-sided ideals is not sufficient: 
\begin{Ex}
Let $R$ be a finite local ring that is not QF, i.e., $S_r$ is not simple as a right ideal, and let $T$ be any simple right ideal. Then $l(T)=J(R)$. Note that because $R$ is local, $top(R)\cong T$ as right $R$-modules. It follows that $|R|/|J(R)|=|T|$, thus proving the Size condition for $T$ and hence for all simple right ideals.     
\end{Ex}

Wood's proof of Thm. \ref{ThmWoodSimple} shows for which simple right ideals $T$ the above inequality holds. In an analogous manner, we can list all simple right ideals $T$ such that the opposite inequality holds, but the method of his proof does not seem to yield this result. We prove this in Prop. \ref{PropQFCardSimple}. First, we recall the following lemma, due to Nakayama, which connects the Nakayama permutation and annihilators. It is taken from the proof of [Nak. 49, Thm. 7], specifically the last paragraph of page 9. 
\begin{Lemma}\label{LemmaAnnofsimple}
Let $R$ be a semiperfect ring with a Nakayama permutation $\pi$ and coinciding socles, and let $J\subseteq I$ be two right ideals such that  $I/J\cong V_k$ for some $k\leq n$.

Then $r(J)/r(I)$ is either zero or a simple module isomorphic to $V'_{\pi^{-1}(k)}$
\end{Lemma}
We omit the proof, as Nakayama's argument remains valid in our more general setting. We inform the reader that the original paper uses $M$ to denote the socle and $\cup$ to denote the sum of two right ideals. 
\begin{remark}
    As the formalism of formal matrix rings gives a description of simple right ideals (Sec. \ref{SecsempSimples}) and left ideals containing the top of the ring [KT 17, Thm. 2.4.3], they can be used to provide an alternative proof of Prop. \ref{PropQFCardSimple}.
\end{remark}
\begin{Prop}\label{PropQFCardSimple}
    Let $R$ be a semiperfect ring with a Nakayama permutation $\pi$ and coinciding socles. 
    
    Then for any $k\leq n$:
 \[
        \begin{cases}
             d(T_k) > d(R/l(T_k)) & \iff \mu_k > \mu_{\pi(k)} \\
             d(T_k) < d(R/l(T_k)) & \iff \mu_k < \mu_{\pi(k)}.
        \end{cases}
    \]
    
    Furthermore, if $R$ is finite, then:
    \[
        \begin{cases}
            |l(T_k)|\cdot|T_k| > |R| & \iff \mu_k > \mu_{\pi(k)} \\
            |l(T_k)|\cdot|T_k| < |R| & \iff \mu_k < \mu_{\pi(k)}.
        \end{cases}
    \]
\end{Prop}
 Given $n<i\leq m$, analogous claim holds for $T_i=soc(e_iR)$ too, using $\mu_{k(i)}$ where $k(i)$ is chosen in such a way that $T_{k(i)}\cong soc(e_iR)$.
\begin{Proof}
We prove the finite case. The claim about generalised dimensions is analogous. 

The right simple ideal $T_k=soc(e_kR)$ is, by assumption, isomorphic to $top(e_{\pi(k)})\cong V_{\pi(k)}$. Then $d(T_k)=\mu_{\pi(k)}$  and hence $|T_k|=|\m_{\pi(k)}|^{\mu_{\pi(k)}}$ (Ex. \ref{ExdimVk}).

Because $l(T_k)$ is not the whole $R$, we have that $l(0)/l(T_k)=R/l(T_k)$ is nonzero and, by the previous lemma, isomorphic to the simple left module $V'_k$. So \[|R|=|l(T_k)||V'_k|=|l(T_k)||\m_k|^{\mu_k}.\]

    Recall that $\m_k\cong \m_{\pi(k)}$ (Prop. \ref{PropIsomfields}), so in particular, they have the same cardinality, which we denote $c$. In total \[|l(T_k)||T_k|=\frac{|R|}{|\m_k|^{\mu_k}} |\m_{\pi(k)}|^{\mu_{\pi(k)}}=|R|c^{(\mu_{\pi(k)}-\mu_k)}. \]
\end{Proof}
\begin{Ex}
    Let $S$ be a local artin algebra and $E$ its self-duality context. Consider trivial matrix rings 
\begin{gather*}
    B:=\begin{pmatrix}
    S&E&0\\
    0&S&E\\
E&0&S \end{pmatrix}\qquad R:=\begin{pmatrix}
    S&E&0&0\\
    0&S&E&E\\
E&0&S&S\\
E&0&S&S
\end{pmatrix}.
\end{gather*}  

Then $B$ is a basic Frobenius ring with a Nakayama permutation $\pi=(1~2~3)$ and $R$ is a ring Morita equivalent to it with multiplicities $\mu_1=\mu_2=1$ and $\mu_3=2$, and associated index sets $I_1=\{1\}, I_2=\{2\}$ and $I_3=\{3,4\}$. 

In the ring $R$ we then have
\begin{gather*}
    l(T_1)=\begin{pmatrix}
    J(S)&E&0&0\\
    0&S&E&E\\
E&0&S&S\\
E&0&S&S
\end{pmatrix} \qquad T_1=\begin{pmatrix}
        0&U_2&0&0\\
        0&0&0&0\\
        0&0&0&0\\
        0&0&0&0\\
    \end{pmatrix}\\[1ex]
    l(T_2)=\begin{pmatrix}
    S&E&0&0\\
    0&J(S)&E&E\\
E&0&S&S\\
E&0&S&S
\end{pmatrix} \qquad T_2=\begin{pmatrix}
        0&0&0&0\\
        0&0&U_3&U_3\\
        0&0&0&0\\
        0&0&0&0\\
    \end{pmatrix}\\[1ex] %
    l(T_3)=\begin{pmatrix}
    S&E&0&0\\
    0&S&E&E\\
E&0&J(S)&S\\
E&0&J(S)&S
\end{pmatrix} \qquad T_3=\begin{pmatrix}
        0&0&0&0\\
        0&0&0&0\\
        U_1&0&0&0\\
        0&0&0&0\\
    \end{pmatrix}\\[1ex]
    l(T_4)=\begin{pmatrix}
    S&E&0&0\\
    0&S&E&E\\
E&0&S&J(S)\\
E&0&S&J(S)
\end{pmatrix} \qquad T_4=\begin{pmatrix}
        0&0&0&0\\
        0&0&0&0\\
        0&0&0&0\\
        U_1&0&0&0\\
    \end{pmatrix}
\end{gather*}
\end{Ex}
Combining the previous proposition with Theorem \ref{ThmAnn1}, we obtain:
\begin{Thm}\label{ThmAnnMain}
    Let $R$ be a semiperfect ring with a Nakayama permutation $\pi$ such that $S_l=S_r$. TFAE

    (1) $\pi$ respects the multiplicities of $R$

    (2) Semisimple right ideals and semisimple left ideals satisfy the Generalised dimension condition.

    (3) Simple right ideals satisfy the Generalised dimension condition.
\end{Thm}
\begin{Proof}
Recall that for a right semisimple ideal $I$, lattices $\mathcal{L}_R(I)$ and $_R\mathcal{L}(R/l(I))$ are anti-isomorphic (Thm. \ref{ThmAnn1}(4)).

(1)$\implies$(2). Assume that $\pi$ respects the multiplicities and  let $\lambda$ be a composition length of $I$ as witnessed by a composition series
\[   0=I_0\subseteq I_1\subseteq \dots\subseteq  I_{\lambda-1}\subseteq I_\lambda=I.\]
Then we have the following composition series of $R/l(I)$   \[ l(I)\subseteq l(I_{\lambda-1})\subseteq \dots \subseteq l(I_1)\subseteq R.
    \] 
    
    For any $1\leq i\leq \lambda$ there is $k(i)$ such that $V_{k(i)}=I_i/I_{i-1}$ and by Lemma \ref{LemmaAnnofsimple}, there is isomorphism of left modules $V'_{\pi^{-1}(k(i))}\cong l(I_{i-1)})/l(I_i)$. 

    This is a bijection between the composition factors of $R/l(I)$ and the composition factors of $I$. Because $\pi$ respects multiplicities, each composition factor is mapped to a simple module of the same generalised dimension. Thus $I$ and  $R/l(I)$ have the same generalised dimension. The proof for left semisimple ideals is analogous.

(2)$\implies$(3) is immediate.

(3)$\implies$(1) Let $T$ be a simple right ideal isomorphic to $V_k$ for some $k\leq n$. Now $R/l(T)$ is a simple left ideal (Lemma \ref{LemmaBijectionSimple}) is isomorphic to $V'_{\pi^{-1}(k)}$ (Lemma \ref{LemmaAnnofsimple}). It follows from the Generalised dimension condition for $T$ that $\mu_k=d(V_k)=d(V'_{\pi^{-1}(k)})=d(V_{\pi^{-1}(k)})=\mu_{\pi^{-1}(k)}$ (Ex. \ref{ExdimVk}). 
    \end{Proof}

\subsection{Results of Honold - Size condition for a general semiperfect ring } \label{SecCardHon}

This section examines cardinality conditions without assuming the Double annihilator property. Lamprecht proved that the existence of a generating character implies the Size condition [Lam. 53, Sec. 3.3 Cor. 1]. Using Wood's result that finite Frobenius rings are characterised by having a generating character [Wood 99, Thm. 3.10],  Honold then proved the converse (Thm. \ref{ThmHonCard}). We reformulate and prove Honold's results for semiperfect rings with essential socles (Cor. \ref{CorCardmain}). In this more general setting, we cannot use the Generalised dimension condition on the Jacobson radical (as a left ideal) and maximal right ideals, as $R$ factored by their respective annihilators does not have a finite length. Instead, we use the right Kasch property and the Generalised dimension condition on the right socle (as a left ideal).

\smallskip

Recall that two semisimple modules are isomorphic if and only if (any) decompositions of them as direct sums of simple submodules have all same multiplicities for all isomorphism types of simple modules. The following Lemma then follows directly from the Definition \ref{DefMult}.
\begin{Lemma}\label{LemmaNakMult}
    Let $R$ be a ring with a Nakayama permutation respecting multiplicities. 

    Then $S_r\cong top(R_R)$.
\end{Lemma}
\begin{remark}
    If R is an artin algebra, then it is a Frobenius ring (Thm. \ref{ThmIovCyclic}).
\end{remark}

This isomorphism is also implied by the Generalised dimension condition:

\begin{Lemma} \label{LemmaHon}Let $R$ be a right Kasch semiperfect ring. Further suppose its right socle satisfies the Generalised dimension condition as a left ideal, and its homogeneous components satisfy it as right ideals. 

    Then $S_r\cong top(R_R)$.
\end{Lemma}
Recall that $S_k$ denotes the homogeneous component of $S_r$ consisting of copies of the simple right module $V_k$. Let us denote its composition length by $n_k$.
\begin{Proof}
     Each $S_k$ is nonzero by the right Kasch assumption and is a right module over the corresponding component $M_{\mu_k}(\m_k)$ of $top(R)$ and hence we represent it as the group of matrices $M_{\mu_k\times n_k}(\m_k)$. In particular $End_R(S_k)\cong M_{n_k}(\m_k)$ as rings. Hence $d(End(S_k))=n_k^2$, whereas $d(S_k)=n_k \mu_{k}$.

    Now for any right ideal, and $S_k$ in particular, we get a homomorphism $\phi\colon R_R\to Hom(S_k, R)$ mapping $r$ to left multiplication by $r$. Because $R$ is Kasch, then $Hom(S_k, R)$ is nonzero. The kernel of this morphism is $l(S_k)$ and hence $\phi(R)\cong R/l(S_k)$. Then we get 
    \begin{align*}
    d(S_k)=d(R/l(S_k))=d(\phi(R))\leq d(Hom(S_k, R)). \\
\intertext{Since the image of $S_k$ under any homomorphism into $R$ is contained in  $S_k$, we get}
        n_k\mu_k= d(S_k)\leq d(Hom_R(S_k, R)) =d(End_R(S_k)) =n_k^
    2.
\end{align*}
$R$ is right Kasch, so $n_k\neq 0$, showing $\mu_k\leq n_k$.

    Because the ring is Kasch, then $r(S_r)=J(R)$ (Lemma \ref{LemmaAnnSoc}). So the Generalised dimension condition implies  \[\sum n_k\mu_k=d(S_r)=d(R/l(S_r))=d(top(R))=\sum \mu_k^2. \] But because $\mu_k \leq n_k$ for each $k$, we get that the only way these two sums can be equal is that  $\mu_k=n_k$ for each $k\leq n$.
\end{Proof}

If $R$ is not an artin algebra, the converse of Lemma \ref{LemmaNakMult} might fail even if $R$ is artinian (Sec. \ref{SecCardClass}). We have, however, the following generalisation of Prop. \ref{PropLamQF}.
\begin{Prop}\label{PropHon}
Let $R$ be a semiperfect ring with an essential right socle and assume that all left primitive ideals have a nonzero socle.

If $S_l\cong top(_RR)$ and $S_r\cong top(R_R)$, then $R$ has a Nakayama permutation respecting multiplicities.     
\end{Prop}
\begin{Proof}
    The isomorphisms imply, in particular, that the ring is Kasch. To show that it has a Nakayama permutation, it is enough to prove that it is a QF-2 ring (Thm. \ref{PropQFKaschEss}). We show that it is the right QF-2; the left version is analogous.
    
    By assumption, we have the following isomorphism of right $R$-modules
\[\bigoplus_{i=1}^mV_i=\bigoplus_{i=k}^k (V_k)^{\mu_k}\cong top(R_R)\cong S_r\cong \bigoplus_{k=1}^n S_k\cong \bigoplus_{k=1}^n (soc(e_kR))^{\mu_{k}}.\]
 The composition length of the left-hand side is $m=\sum \mu_k$. Because $S_r$ is essential, we see that each $soc(e_kR)$ is nonzero. Therefore, the composition length of the right-hand side is at least $m$, forcing it to be exactly $m$. Hence, all primitive right ideals have simple socles, which proves that $R$ is a right QF-2 ring.  

Let us denote the Nakayama permutation of $R$ as $\pi$. Then 
\[\bigoplus_{k=1}^n (V_k)^{\mu_k}\cong  \bigoplus_{k=1}^n (top(e_kR))^{\mu_k} \cong \bigoplus_{k=1}^n (soc(e_{\pi^{-1}(k)}R))^{\mu_k} \]
and because $V_k$ are pairwise nonisomorphic, comparing the composition lengths of the left-hand side and right-hand side forces $\mu_k=\mu_{\pi^{-1}(k)}$.
\end{Proof}
In total, we get:
\begin{Cor}\label{CorCardmain}
Let $R$ be a semiperfect ring with essential socles. The following are equivalent:

(1) $R$ has a Nakayama permutation respecting the multiplicities.

(2) $S_l\cong top(_RR)$ and $S_r\cong top(R_R)$.

(3) $R$ is right Kasch, and its right (left) socle satisfies the Generalised dimension condition as a left (right) ideal, and its homogeneous components satisfy it as right (left) ideals.
\end{Cor}
\begin{Proof}
(3)$\implies$(2) follows from Lemma \ref{LemmaHon}. Implication (2)$\implies$(1) follows from Prop. \ref{PropHon}. Assuming (1), $R$ is Kasch by definition. Furthermore, $S_l=S_r$ (Thm. \ref{PropQFKaschEss}), so the rest follows by  Thm.  \ref{ThmAnnMain}.
\end{Proof}

\textbf{Acknowledgments:} 
The author would like to thank Steven Dougherty for discussions during NCRA IX that provided the impetus to compile these results into the present survey format."

\bigskip

\noindent \textbf{Bibliography}

[AF 92] F.  W. Anderson, K. R. Fuller. (1992). Rings and Categories of Modules. Graduate Texts in Mathematics 13. (Springer). https://doi.org/10.1007/978-1-4612-4418-9

[Art. 27] E. Artin (1927). Zur Theorie der hyperkomplexen Zahlen. \textit{Abh. Math. Semin. Univ. Hambg}. 5. 251–260. https://doi.org/10.1007/BF02952526 

[AW 09] D. D. Anderson, M. Winders. (2009). Idealization of a Module. \textit{J. Com. Alg.} 1(1). 3-56. https://doi.org/10.1216/JCA-2009-1-1-3

[Azu. 59] G. Azumaya. (1959). A Duality Theory for Injective Modules. (Theory of Quasi-Frobenius Modules). \textit{Am. J. Math.} 81(1). 249-278. https://doi.org/10.2307/2372855

[Bir. 48] G. Birkhoff (1948). Lattice Theory, Revised Edition. (American Mathematical Society). ISBN: 0821889532.

[Bass 60] H. Bass. (1960). Finitistic Dimension and a Homological Generalization of Semi-Primary Rings. \textit{Trans. Am. Math. Soc.} 95(3). 466-488. https://doi.org/10.2307/1993568

[Baer 43] R. Baer. (1943). Rings with Duals. \textit{Am. J. Math}. 65(4). 569–584. https://doi.org/10.2307/2371866

[Baer 52] R. Baer. (1952). Linear Algebra
and Projective Geometry. Pure and Aplied Mathematics 2. (Academic Press) 

[BO 09] Y. Baba, K. Oshiro. (2009). Classical Artinian Rings and Related Topics. (World Scientific Publishing). https://doi.org/10.1142/7451

[CDY 04] J. Chen, N. Ding \& M. Yousif. (2004). On Noetherian Rings with Essential Socle. \textit{J. Aust. Math. Soc.} 76(1). 39-50. https://doi.org/10.1017/S1446788700008685

[CG 92] H.L. Claasen,  R.W. Goldbach. (1992).  A field-like property of finite rings. \textit{Indag. Math}. 3(1). 11-26. https://doi.org/10.1016/0019-3577(92)90024-F

[Cha. 61] S. U. Chase. (1961). A Generalization of the Ring of Triangular Matrices. \textit{Nagoya Math. J.} 18. 13-25. https://doi.org/10.1017/S0027763000002208

[CNY 00] V. Camillo,  W.K. Nicholson \&  M.F. Yousif (2000). Ikeda–Nakayama Rings. \textit{J. Alg.} 226(2). 1001-1010. https://doi.org/10.1006/jabr.1999.8217

[CR 64] C. W. Curtis, I. Reiner. (1964). Representation Theory of Finite Groups and Associative Algebras. (Interscience). https://doi.org/10.1017/S0008439500026990

[Die. 58] J. Dieudonné. (1958). Remarks on quasi-Frobenius rings. Ill. J. Math. 2(3). 346-354. https://doi.org/10.1215/ijm/1255454538

[Dou. 17] S. T. Dougherty.  (2017). Algebraic Coding Theory Over Finite Commutative Rings. SpringerBriefs in Mathematics. (Springer). https://doi.org/10.1007/978-3-319-59806-2.

[Fai. 66] C. Faith. (1966). Rings with Ascending Condition on Annihilators.\textit{ Nagoya Math. J.} 27(1). 179-191. https://doi.org/10.1017/S0027763000011983

[Fai. 76] C. Faith. (1976). Algebra II:
Ring Theory. Grundlehren der mathematischen Wissenschaften 191. (Springer). https://doi.org/10.1007/978-3-642-65321-6  

[Ful. 69] K. Fuller. (1969). On indecomposable Injectives over Artinian Rings. \textit{Pac. J. Math.} 29(1). 115-135 https://doi.org/10.2140/pjm.1969.29.115

[FV 02] C. Faith, D. Van Huynh. (2002).  When Self-injective Rings are QF: a Report on a Problem. 	\textit{J. Algebra Its Appl.} 1(1). 75-105. https://doi.org/10.1142/S0219498802000070

[FW 67]  C. Faith, E. A. Walker. (1967). Direct-Sum Representations of Injective Modules. J. Alg. 5(2). https://doi.org/10.1016/0021-8693(67)90035-X

[Gre. 82] E. L. Green. (1982). On the Representation Theory of Rings in Matrix Form.  \textit{Pac. J. Math.} 100(1). 123-138. https://doi.org/10.2140/pjm.1982.100.123

[GNW 04] M. Greferath, A. Nechaev \& R. Wisbauer. (2004). Finite Quasi-Frobenius Modules and Linear Codes. \textit{J. Algebra Appl.} 3(3). 247-272. https://doi.org/10.1142/S0219498804000873.

[Har. 82] M. Harada. (1982). Self mini-injective rings. \textit{Osaka J. Math.} 19(3). 587-597.

[HN 85] C. R. Hajarnavis, N. C. Norton. (1985). On Dual Rings and their Modules. \textit{J. Alg.} 93(2). 253-266. https://doi.org/10.1016/0021-8693(85)90159-0

[Han. 73a] T. A. Hannula. (1973). On the Construction of Quasi-Frobenius Rings. \textit{J Alg.} 25(3). 403-413. https://doi.org/10.1016/0021-8693(73)90089-6

[Han. 73b]  T. A. Hannula. (1973). The morita context and the construction of QF rings. In: \textit{Proceedings of the Conference on Orders, Group Rings and Related Topics}. Lecture Notes in Mathematics 353. (Springer). https://doi.org/10.1007/BFb0059264

[Hir. 97] Y. Hirano. (1997). On admissible rings. \textit{Indag. Math}. 8(1).  55-59. https://doi.org/10.1016/S0019-3577(97)83350-2

[Hon. 01] T.  Honold. (2001). Characterization of finite Frobenius rings. Arch. Math. 76(6). 406-415. https://doi.org/10.1007/PL00000451

[Hop. 38] C. Hopkins. (1938). Nil-rings with minimal condition for admissible left ideals. \textit{Duke Math. J.} 4(4). 664-667. https://doi.org/10.1215/S0012-7094-38-00457-0

[Hop. 39] C. Hopkins. (1939). Rings With Minimal Condition for Left Ideals. \textit{Ann. Math.} 40(3). 712-730. https://doi.org/10.2307/1968951

[Ike. 51] M. Ikeda. (1951). Some generalizations of quasi-Frobenius rings. \textit{Osaka Math. J.} 3(2). 227-239. https://ir.library.osaka-u.ac.jp/repo/ouka/all/10382/omj03\_02\_05.pdf 

[Ike. 52] M. Ikeda. (1952). A characterisation of quasi-Frobenius rings. \textit{Osaka Math. J.} 4(2). 203-209.

[IN 54] M. Ikeda, T. Nakayama. (1954). On Some Characteristic Properties of Quasi-Frobenius and Regular Rings. \textit{Proc. Am. Math. Soc.} 5(1). 15-19. https://doi.org/10.2307/2032097

 [Iov. 16] M. C. Iovanov. (2016). Frobenius–Artin algebras and infinite linear codes. \textit{J. Pure Appl. Algebra}. 560-576. 
 220(1). https://doi.org/10.1016/j.jpaa.2015.05.030

  [Iov. 22] Iovanov, M. C. (2022). On Infinite MacWilliams Rings and Minimal Injectivity Conditions. \textit{Proc. Am. Math. Soc.} 150(11). 4575-4586. https://doi.org/10.1090/proc/15929

[Kato 67] T. Kato. (1967). Self-injective rings. \textit{Tahoku Math. J.} 19(4).  485-495. https://doi.org/10.2748/tmj/1178243253

[Kas. 54] F. Kasch. (1954). Grundlagen einer Theorie der Frobeniuserweiterungen. \textit{Math. Ann.} 127. 453–474. https://doi.org/10.1007/BF01361137

[Kas.  82] F. Kasch. (1982). Modules and rings.  LMS monograph 17. (Academic Press). https://doi.org/10.1017/S001309150001703X.

[Kra. 24] D. Krasula. Möbius function for modules and thin representations. \textit{ArXiv}. https://doi.org/10.48550/arXiv.2403.05656

[Kra. 25] D. Krasula. (2025). Formal matrix representations of rings with a Nakayama permutation. \textit{ArXiv}

[Kra. 26] D. Krasula. (2026). Endomorphism rings of Simple Modules and Block Decomposition. \textit{J. Algebra Its Appl.}  Online Ready https://doi.org/10.1142/S0219498826501562

[Krs. 07] H. Krause. (2007). An axiomatic characterization of the Gabriel-Roiter measure.\textit{ Bull. London Math.} 39(4). 550-558. https://doi.org/10.1112/blms/bdm033

[KT 17] P. Krylov, A. Tuganbaev. (2017). Formal Matrices. Algebra and Applications 23. (Springer). https://doi.org/10.1007/978-3-319-53907-2

[Lam 91] T. Y. Lam. (1991). A First Course in Noncommutative Rings.  Graduate Texts in Mathematics 131. (Springer). https://doi.org/10.1007/978-1-4419-8616-0

[Lam 99] T. Y. Lam. (1999). Lectures on Modules and Rings.  Graduate Texts in Mathematics 189. (Springer). https://doi.org/10.1007/978-1-4612-0525-8

[Lam. 53] E. Lamprecht. (1953). Über I-reguläre Ringe, reguläre Ideale und Erklärungsmoduln. I. Math. Nachr. 10(5-6). 353-382. https://doi.org/10.1002/mana.19530100506

 [Lev. 40] J. Levitzki. (1940). On rings which satisfy the minimum condition for the right-hand ideals. \textit{Compos. Math.} 7. 214–222.

[McD. 74] B. R. McDonald. (1974). Finite Rings with Identity. Pure and applied mathematics (M. Dekker).

[Mor. 66] K. Morita. (1966). On S-Rings in the Sense of F. Kasch. \textit{Nagoya Math. J.} 27(2). 687-695.  https://doi.org/10.1017/S0027763000026477

[Mül. 70] B. J. Müller. (1970). On Semi-perfect Rings.  \textit{Illinois J. Math} 14(3). 464-467. https://doi.org/10.1215/ijm/1256053082

[Mül. 74] B. J. Müller. (1974). The Structure of Quasi-Frobenius Rings. \textit{Can. J. Math.} 26(5). 1141-1151. https://doi.org/10.4153/CJM-1974-106-7

[Nak. 39] T. Nakayama. (1939). On Frobeniusean Algebras. I. \textit{Ann. Math.} 40(3). 611-633. https://doi.org/10.2307/1968946

[Nak. 41a]  T. Nakayama. (1941). On Frobeniusean Algebras II. \textit{Ann. Math.} 42(1). 1-21. https://doi.org/10.2307/1968984

[Nak. 41b] T. Nakayama. (1941). 
Algebras with anti-isomorphic left and right ideal lattices. \textit{Proc. Imp. Acad.} 17(3). 53-56  https://doi.org/10.3792/pia/1195578880

[Nak. 49] T. Nakayama. (1949). Supplementary remarks on Frobeniusean algebras I. \textit{Proc. Japan Acad.} 25(7). 45-50. https://doi.org/10.3792/pja/1195571908

[NKC 23] T. H. Nhan,  M.  Tamer Koşan \&  T. Cong Quynh (2023). On Annihilators and Quasi-Frobenius Rings. \textit{Lobachevskii J. Math.} 44. 5355–5363. https://doi.org/10.1134/S1995080223120259

[Nor. 75] N. C. Norton. (1975). Generalizations of the theory of quasi-Frobenius rings.  \textit{Ph.D. thesis}. University of Warwick

[NY 97] W. Nicholson, M. Yousif. (1997). Mininjective Rings. \textit{J. Alg.} 187(2). 548-578. https://doi.org/10.1006/jabr.1996.6796

[NY 03] W. Nicholson, M. Yousif. (2003). Quasi-Frobenius Rings. 
Cambridge Tracts in Mathematics 158. (Cambridge University Press).
https://doi.org/10.1017/CBO9780511546525

[Roux 71] B. Roux. (1971). Sur la Dualité de Morita. \textit{Tahoku Math. J.} 23(3). 457-472. https://doi.org/10.2748/tmj/1178242594

[Sto. 69] H. H. Storrer. (1969). \textit{Can. Math. Bull.} 12(3). 287-292.  https://doi.org/10.4153/CMB-1969-036-9

[Thr. 48] R. M. Thrall. (1948). Some Generalizations of Quasi-Frobenius Algebras. \textit{Trans. Am. Math Soc.} 64(1). 173-183. https://doi.org/10.2307/1990561

[Wood 99] J. A. Wood. (1999). Duality for Modules over Finite Rings and Applications to Coding Theory. \textit{ Am. J. Math.} 121(3). 555-575. https://doi.org/10.1353/ajm.1999.0024

[Wood 09]  J. A. Wood. (2009). Foundations of Linear Codes Defined over Finite Modules: the Extension Theorem and the MacWilliams Identities. In: Codes Over Rings. https://doi.org/10.1142/9789812837691\_0004

[Xue 89] W. Xue. (1989). Two Examples of Local Artinian Rings. \textit{Proc. Am. Math. Soc.} 107(1). 63-65. https://doi.org/10.2307/2048035

[Xue 92] W. Xue. (1992). Rings with Morita duality. Lecture Notes in Mathematics 1523. (Springer). https://doi.org/10.1007/BFb0089071

[Xue 97]  W. Xue. (1997). On a Theorem of Fuller. \textit{J. Pure Appl. Algebra}. 122(1-2). 159-168. https://doi.org/10.1016/S0022-4049(96)00070-9

[Xue 98] W. Xue. (1998). A note on finite local rings. \textit{Indag. Math.} 9(4). 627-628. https://doi.org/10.1016/S0019-3577(98)80040-2

\end{document}